\documentclass[a4paper,11pt,leqno]{article}
\usepackage{theorem}
\usepackage{latexsym}
\usepackage{amsmath}
\usepackage{delarray}
\usepackage{graphicx}
\usepackage{delarray}
        \usepackage[latin1]{inputenc}
        \usepackage{color}
        \usepackage{array}
        \usepackage{longtable}
        \usepackage{calc}
        \usepackage{multirow}
        \usepackage{hhline}
        \usepackage{ifthen}

\makeatletter
\@addtoreset{equation}{section}
\makeatother

\theorembodyfont{\rmfamily}

\def\sgn{\mathop{\rm sgn}}

\let\margin\marginpar
\newcommand\myMargin[1]{\margin{
      \raggedright\tiny  \setlength{\baselineskip}{0pt} #1}}
\renewcommand{\marginpar}[1]{\myMargin{#1}}

\title{Some improvements in numerical evaluation of symmetric stable density and
  its derivatives}
\author{Muneya MATSUI${}^{\dag}$ and Akimichi TAKEMURA${}^{\ddag 
}$ \\ \\
$\dag$ Graduate School of Economics, University of Tokyo
 \\
$\ddag$ Graduate School of Information Science and Technology,
University of Tokyo
}
\date{}
\begin{document}
\maketitle

\begin{abstract}
  We propose improvements in numerical evaluation of symmetric stable
  density and its partial derivatives with respect to the parameters.
  They are useful for more reliable evaluation of maximum likelihood
  estimator and its standard error.  Numerical values of the Fisher
  information matrix of symmetric stable distributions are also given.
  Our improvements consist of modification of the method of Nolan
  (1997) for the boundary cases, i.e., 
  in the tail and mode of the densities and in the neighborhood of the
  Cauchy and the normal distributions.
\end{abstract}

\section{Introduction}
There have been many researches to evaluate densities and quantiles
of symmetric or general stable distributions. 
McCulloch (1998) considered efficient algorithms for
approximating symmetric stable densities $f(x;\alpha)$ for the range
$\alpha>0.85$, where parameter $\alpha$ denotes the characteristic exponent.
Nolan (1997) gave accurate algorithms for
general stable densities based on 
integral representations of the densities which were derived by 
Zolotarev (1986).  
Nolan provides a very useful program package ``STABLE'' 
on his web page\footnote{{\tt http://academic2.american.edu/\~{}jpnolan}}.
However his program exhibits some unreliable behavior around the boundary  
as stated in the users
guide of STABLE.
Therefore even in the case of symmetric stable distributions, reliable
computations of density functions including all the boundary cases is 
still needed. Furthermore for maximum likelihood estimation, it is
desirable to directly compute the derivatives of the density function.
In this paper we present reliable computations of symmetric stable
density functions and their partial derivatives.  Our computation of
densities is accurate for all values of $x$ and $0.1 < \alpha \le 2$.
Concerning the partial derivatives it is accurate in a somewhat smaller
range of values.

Regarding maximum likelihood estimation for the range $\alpha \geq 0.4$,
Nolan (2001) used interpolated stable
densities and maximized the likelihood by approximate
gradient search (constrained quasi-Newton method) because of its
efficiency. But near the boundary of the parameter space 
interpolation may be inaccurate and the direct integral representation
is  used.
Note that the direct integral representation is also not very reliable
and slow near the boundary. 
In the symmetric case, Brorsen and Yang (1990) discussed maximum likelihood  estimation
using an integral representation of the densities given by Zolotarev (1986).
But they have only considered the range $\alpha>1$ to avoid the
discontinuity and nondifferentiability at $\alpha=1$. 
Furthermore they did not check the sample covariances of their maximum
likelihood computation with the Fisher information matrix.

These previous researches on maximum likelihood estimation 
have not used direct evaluation of the derivatives of the log
likelihood function with respect to the parameters. 
For reliable evaluation of the maximum likelihood estimator and its
standard error, direct and reliable evaluation of the first and the second
derivatives of the log likelihood function is desirable.

The organization of this paper is as follows.
In Section \ref{sed:preliminaris} we summarize notations and
preliminary results on symmetric stable density.
In Section \ref{sec:density} we provide an accurate algorithm for calculations 
of symmetric stable distributions
which modifies Nolan (1997) for $x$ near $0$ or $\infty$ and
for $\alpha=1$ or $\alpha=2$ using various expansions.
Accurate algorithms for the partial derivatives of symmetric stable
distributions with respect to the parameters are given in Section 
\ref{sec:derivative-location-scale}
and
Section \ref{sec:derivative-alpha}. 
Fisher information matrices are calculated in Section \ref{sec:fisher-information},
together with some simulation studies on the variance of the maximum
likelihood estimator and the observed Fisher information. We also
discuss behavior of Fisher information as $\alpha \rightarrow 2$.
Some discussions are given in Section \ref{sec:discussion}.  

\section{Preliminaries}
\label{sed:preliminaris}

In this section we prepare notations and summarize preliminary
results. There are many parameterizations for stable distributions
and much confusion has been caused. Our parameterizations follow the useful 
parameterizations for statistical inference which were given in Nolan (1998).

\subsection{Notations}
\label{sec:notation}
Let
\begin{equation}
\label{eq:characteristic-function}
\Phi(t)=\Phi(t;\mu,\sigma,\alpha)=\exp\bigl(-|\sigma t|^{\alpha} + i
\mu t\bigl)
\end{equation}
denote the characteristic function of symmetric stable distribution
with parameters
$$
\theta=(\theta_1, \theta_2, \theta_3)=(\mu,\sigma,\alpha),
$$
where $\alpha$ is the characteristic exponent, $\mu$ is a location
parameter and $\sigma$ is a scale parameter.
For the standard case
$(\mu,\sigma)=(0,1)$ we simply write the characteristic function as
$
\Phi(t;\alpha)=\exp(-|t|^{\alpha}).$

The corresponding density is written as
$f(x;\mu,\sigma,\alpha)$  and $f(x;\alpha)$ in the standard case:
$$
f(x;\mu,\sigma,\alpha) =
\frac{1}{\sigma}f(\frac{x-\mu}{\sigma};\alpha).
$$
At $\alpha=1$ 
\[
f(x;1)= %
   \frac{1}{\pi(1+x^2)}
\]
is the Cauchy density and  at $\alpha=2$ 
\[
f(x;2)= %
  \frac{1}{2\sqrt{\pi}}\exp(-x^2/4)
\]
is the normal density $N(0,2)$.
Accordingly the density can be defined to 
constitute a location-scale family.
In the following, the first derivative of $f(x;\alpha)$ with respect to $x$ 
is denoted by $f'(x;\alpha)$ and the second derivative is denoted by
$f''(x;\alpha)$. 
Then 
\begin{eqnarray}
\label{eq:f-mu}
&&
\frac{\partial}{\partial
  x}f(x;\mu,\sigma,\alpha)=-\frac{\partial}{\partial
  \mu}f(x;\mu,\sigma,\alpha)=\frac{1}{\sigma^2} f'(\frac{x-\mu}{\sigma};\alpha),\\
&&
\frac{\partial^2}{\partial
  x^2}f(x;\mu,\sigma,\alpha)=\frac{\partial^2}{\partial
  \mu^2}f(x;\mu,\sigma,\alpha)=\frac{1}{\sigma^3} f''(\frac{x-\mu}{\sigma};\alpha).
\end{eqnarray}
The partial derivatives with respect to $\sigma$ and $\alpha$ are
written by subscripts, e.g., 
$$
f_\alpha(x;\mu,\sigma,\alpha)
=\frac{\partial}{\partial \alpha} f(x;\mu,\sigma,\alpha), 
\quad f_{\alpha\sigma} (x;\mu,\sigma,\alpha)=
\frac{\partial^2}{\partial\alpha\partial\sigma}
f(x;\mu,\sigma,\alpha).
$$
As above, when these derivatives are evaluated at the standard case
$(\mu,\sigma)=(0,1)$ we write 
$f_\alpha(x;\alpha)$, $f_{\alpha\sigma} (x;\alpha)$, etc.  Note that
\begin{equation}
\label{eq:f-sigma}
f_\sigma(x;\alpha)= -f(x;\alpha) - x f'(x;\alpha),
\quad 
f_{\sigma\sigma}(x;\alpha)= 2f(x;\alpha) +4x f'(x;\alpha)+ x^2f''(x;\alpha).
\end{equation}
Furthermore we write %
$$
f_{\alpha}'(x ; \alpha)=\frac{\partial^2}{\partial x
  \partial\alpha}f(x ; \alpha).
$$

The reason we consider up to the second order derivatives of the density
function is that in assessing the standard error of the maximum
likelihood estimator $\hat \theta$, the observed Fisher information
\begin{equation}
\label{eq:observed-FI}
\hat I_{\theta\theta}(x_1,\ldots,x_n) = - \frac{1}{n}\sum_{i=1}^n \frac{\partial^2 \log
  f(x_i ; \hat\theta)}{\partial \theta^2}
= \frac{1}{n}\sum_{i=1}^n \left(\frac{ \frac{\partial}{\partial\theta}
  f(x_i ; \hat\theta)}{f(x_i ; \hat\theta)}\right)^2
 - \frac{1}{n}\sum_{i=1}^n \frac{ \frac{\partial^2}{\partial\theta^2}
  f(x_i ; \hat\theta)}{f(x_i ; \hat\theta)}
\end{equation}
is usually preferred to the value of the Fisher information matrix at
$\hat\theta$ (e.g. Efron and Hinkley (1978)).

Note that there are other parameterizations of symmetric stable
distributions than (\ref{eq:characteristic-function}).  However
different parameterizations in the literature are smooth functions of
each other including the boundary $\alpha=2$ and differentiations in
terms of other parameterizations can be obtained from the results of
this paper by the chain rule of differentiation.

\subsection{Preliminary results}

{}From  equation (2.2.18) of %
Zolotarev (1986) %
or Theorem 1 of Nolan (1997), 
the density $f(x;\alpha)$ for the case 
$\alpha\neq 1$ and $x>0$  is written as
\begin{equation}
\label{dense}
f(x;\alpha)
 = \frac{\alpha}{\pi|\alpha-1|x}\int^{\frac{\pi}{2}}_0 
     g(\varphi;\alpha,x)\exp(-g(\varphi;\alpha,x)) d\varphi, 
\end{equation}
where 
\begin{equation}
g(\varphi;\alpha,x)=\left(\frac{x\cos \varphi}{\sin \alpha
    \varphi}\right)^{\frac{\alpha}{\alpha-1}}\frac{\cos(\alpha-1)\varphi}{\cos \varphi}. 
\label{integrand1}
\end{equation}
Note that at  $x=0$
\[
f(0;\alpha)= \frac{1}{\pi}\Gamma\left(1+\frac{1}{\alpha}\right)
\]
for all $0 < \alpha \le 2$.

For the case $x\to 0$ and  
$0 < \alpha \leq 2$,  $\alpha \neq 1$, the following expansion
can be used.
\begin{equation}
\label{eq:xto0}
f(x;\alpha)=
\frac{1}{\pi}\sum^{\infty}_{k=0}\frac{\Gamma((2k+1)/\alpha+1)}{(2k+1)!}(-1)^k
x^{2k} 
=
\frac{1}{\pi\alpha}\sum^{\infty}_{k=0}\frac{\Gamma((2k+1)/\alpha)}{(2k)!}(-1)^k
x^{2k}. 
\end{equation}
For  $\alpha < 1$, this series is not convergent but can be
justified as an asymptotic expansion as
$x\rightarrow 0$.  For $1 < \alpha \le 2$, it 
is convergent for every $x$.
Similarly for 
the case $x\to \infty$ and  
$0 < \alpha \le 2$,  $\alpha \neq 1$, we have
\begin{equation}
 \label{eq:xtoinf}
f(x;\alpha)=
  \frac{1}{\pi}\sum^{\infty}_{k=1}\frac{\Gamma(k\alpha+1)}{k!} (-1)^{k-1}
 \sin(\frac{\pi \alpha k }{2}) x^{-k\alpha-1}. 
\end{equation}
For $\alpha < 1$ this series converges for every $x\neq 0$ and 
for $\alpha>1$ this series can be justified as an asymptotic expansion
as
$x\rightarrow\infty$. For $\alpha=2$ this asymptotic expansion 
is zero, which corresponds to the fact that the tail of
normal distribution is exponentially small.
These  (asymptotic) expansions are stated in Bergstr\"om
(1953), Section XVII.6 of Feller (1971), Section 2.4 and Section 2.5 of
Zolotarev (1986).

\section{Numerical evaluation of symmetric stable densities}
\label{sec:density}

As in Nolan (1997) we numerically evaluate the density function using (\ref{dense}).
In (\ref{dense}) 
the function $g(\varphi;\alpha,x):
[0,\frac{\pi}{2}] \to [0,\infty]$ plays an important role, because  
properties of this function make the numerical integration quite efficient. Note that 
$g(\varphi;\alpha,x)$ is continuous and positive, strictly increases
from $0$ to $\infty$ for $\alpha<1$ and strictly decreases from
$\infty$ to $0$ for $\alpha>1$. Therefore the integrand
$g(\cdot)\exp(-g(\cdot))$ is unimodal and its maximum value $1/e$ is
uniquely attained  at $\varphi_1$ satisfying
$g(\varphi_1;\alpha,x)=1$.  
When the value of
density is small, the integrand  concentrates around its mode very narrowly.
Then quadrature algorithms  may miss the
integrand.  Therefore we solve $g(\varphi;\alpha,x)=1$ for $\varphi_1$
and the integral is divided into two intervals around this mode (see
Nolan (1997)).  
For numerical calculations of
(\ref{dense}) we use adaptive integration with singularity (QAGS) in
GNU Scientific Library (2003).  %
For most values of $x$
and $\alpha$ this integration works well.

However, when $\alpha$ is close to $2$ this algorithm has some difficulty. 
Note that for $\alpha=2$
\[
 g(\varphi;2,x)=\left( \frac{x }{2 \sin \varphi} \right)^2
\]
and $g(\pi/2;2,x)=x/2$. Therefore 
$\varphi_1$ exceed $\pi/2$ when $x>2$.
There are some other numerical difficulties in (\ref{dense}). 
We list these difficulties and propose alternative practical methods for
evaluating the density. 

\begin{enumerate}
\setlength{\itemsep}{2pt}
\item  $\alpha$ is small and $x \to 0$: \\
  If $\alpha$ is small, the density is very much concentrated at $x=0$.
  For example Nolan (1997) states %
  $f(0;0.1)=1.155 \times 10^6$ whereas $f(0.01;0.1)=1.66$.
  In our calculations when
  $\alpha$ is small and $x$ near $0$ the values of (\ref{dense}) 
  sometimes become larger than $f(0;\alpha)$, contradicting
  the unimodality  of the stable density.  For this case we can use
  the asymptotic expansion (\ref{eq:xto0}).
 
\item $x \to \infty$: \\ 
We cannot guarantee the accuracy of (\ref{dense}) in the case of $x \to \infty$.
Since stable distributions have heavy tails, 
  reliable calculation of their densities is needed for large $x$.
  In our calculations when   $\alpha>1$ and $x$ is large
  the values of (\ref{dense}) sometimes become 
  much smaller than the asymptotic expansion (\ref{eq:xtoinf}). 
  For this case we can use the asymptotic expansion (\ref{eq:xtoinf}).

\item $\alpha$ is near 1: \\
The representation (\ref{dense}) can not be applied at $\alpha=1$
theoretically. The numerical quadrature of (\ref{dense}) becomes
unreliable because of roundoff errors, when $\alpha$ is close to $1$.

\item $\alpha$ is near 2: \\
Though the representation (\ref{dense}) can be applied near $\alpha=2$
theoretically, 
it seems to be too close to the normal distribution 
in the tail of the distribution. 
Actually the values of the density in the tail obtained by the
integral representation (\ref{dense}) %
is much smaller than the asymptotic expansion (\ref{eq:xtoinf}).
\end{enumerate}

For the rest of this section, we discuss the cases 3 and 4 above.

For $\alpha\doteq 1$, we consider Taylor expansion of the density
around $\alpha=1$.
Let 
$$\gamma \doteq 0.57722$$
denote the Euler's constant throughout the rest of this paper.
The Taylor expansion of $f(x;\alpha)$ around 
$\alpha=1$ is given as follows.
\begin{equation}
f(x;\alpha) = f(x;1)+f_{\alpha}(x;1)(\alpha-1) 
+\frac{1}{2}f_{\alpha\alpha}(x;1)(\alpha-1)^2
+\frac{1}{6}f_{\alpha\alpha\alpha}(x;1)(\alpha-1)^3
+ o((\alpha-1)^3)
,\label{eq:alpha-1}
\end{equation}
where
\begin{eqnarray}
f_{\alpha}(x;1) &=& \frac{1}{\pi}
\left\{\frac{x^2-1}{(1+x^2)^2}\left(1-\gamma-\frac{1}{2}\log(1+x^2)\right)
+\frac{2x}{(1+x^2)^2}\arctan x\right\}, 
\label{eq:alpha-1-1}\\
\label{eq:alpha-1-2}
f_{\alpha\alpha}(x;1) &=& 
\frac{x^4-6x^2+1}{\pi(1+x^2)^3}
\left\{\frac{\pi^2}{6}+\left(1-\gamma-\frac{1}{2}\log(1+x^2)\right)^2-1-\arctan^2x\right\} \\
& &+
\frac{8x(x^2-1)}{\pi(1+x^2)^3}\arctan x
\left(\frac{3}{2}-\gamma-\frac{1}{2}\log(1+x^2)\right)
     \nonumber \\
& &+
\frac{2}{\pi(1+x^2)^3}\left\{(1-3x^2)
\left(1-\gamma-\frac{1}{2}\log(1+x^2)\right)-x(1+x^2)\arctan x\right\}.
\nonumber
\end{eqnarray}
For convenience the explicit form of 
$f_{\alpha\alpha\alpha}(x;1)$ is given in appendix
\ref{sec:alpha-3}.
(\ref{eq:alpha-1-1}) and 
(\ref{eq:alpha-1-2}) are proved as follows.
{}From the equation 4.40 on page
18 of Oberhettinger (1990)
\begin{equation}
\label{eq:oberhettinger}
\int_0^\infty  u^{\nu-1}e^{-au} \log u \cos(u y) du =  
(a^2 + y^2)^{-\frac{1}{2}\nu} \Gamma(\nu)\bigl[
  \cos(\nu z) \{\psi(\nu)- \frac{1}{2} \log(a^2 + y^2)\} - z \sin(\nu
  z)\bigl], 
\end{equation}
where 
$$ {\rm Re}\;\nu > 0, \quad z=\arctan(y/a), \quad
\psi(\nu)=\Gamma'(\nu)/\Gamma(\nu).
$$
Differentiating  (\ref{eq:oberhettinger})
several times with respect to $\nu$, setting $\nu=2$, 
and combining the results in the inversion formula we obtain 
(\ref{eq:alpha-1-1}) and 
(\ref{eq:alpha-1-2}).
The conditions for change of integral and differentiation are satisfied in these cases.
Although higher order derivatives of $f(x;\alpha)$ with respect to
$\alpha$ can be evaluated along the same line, we found that the three term
expansion (\ref{eq:alpha-1}) is sufficiently accurate.

Now we consider the case $\alpha \to 2$ and $x \to \infty$. It seems natural to use 
asymptotic expansion (\ref{eq:xtoinf}). However the normal density has
an exponentially small tail and this expansion is meaningless for $\alpha=2$.
However %
in view of smoothness at $\alpha=2$, an 
approximation around the normal density is desirable. 
This case is somewhat subtle, but we found that  the following
procedure works well numerically.
Note that from  (\ref{eq:xto0}),  for each fixed  $x$, $f(x;\alpha)$ is differentiable with
respect to $\alpha$  $(>1)$ even at $\alpha=2$, i.e., for  each fixed $x$ we have
\begin{equation}\label{eq:alpha-2}
f(x;\alpha)= f(x;2)+f_{\alpha}(x;\alpha)(\alpha-2)+\frac{1}{2}
f_{\alpha\alpha}(x;\alpha)(\alpha-2)^2 + o((\alpha-2)^2).
\end{equation}
Differentiating (\ref{eq:xtoinf}) with respect to $\alpha$, for large $x$, heuristically we have
$$ 
 f_{\alpha}(x;\alpha)(\alpha-2) 
\sim 
\frac{1}{\pi}\Gamma(\alpha+1)
\frac{\pi}{2}\cos\left(\frac{\pi\alpha}{2}\right)
x^{-\alpha-1}(\alpha-2) 
\sim 
\frac{\Gamma(\alpha+1)}{2}x^{-\alpha-1}(2-\alpha) 
$$
and 
$$
f(x;\alpha) \sim f_{\alpha}(x;\alpha)(\alpha-2).
$$

We summarize our treatments of various boundary cases in Table
\ref{tdense}.
The range of $x$ and
$\alpha$ and the number of terms $k$ in
the expansions (\ref{eq:xto0}) and (\ref{eq:xtoinf}) are shown.
For $0.99 < \alpha \le 1.01$ we use formula (\ref{eq:alpha-1}) and 
for $\alpha > 1.99999$ and for $x$ large we use  the maximum of
(\ref{dense}) and (\ref{eq:alpha-2}).
Note that for most cases a small number of terms in the expansions
(\ref{eq:xto0}) or (\ref{eq:xtoinf})
is sufficient.
The approximations are very effective since the expansions  and 
the integral (\ref{dense}) give virtually the same results for most values of $x$ and
$\alpha$.

\begin{table}[htbp]
\caption[percentage]{Approximations to stable density at boundary cases} 
\label{tbl:ttdense}
\vspace{-3mm}
\begin{center}
\begin{tabular}{|c|c|c|c|c|} \hline
$\alpha \backslash x$  & \multicolumn{2}{c|}{$x \to 0$, formula   (\ref{eq:xto0})}
   & \multicolumn{2}{c|}{$x \to \infty$,  formula (\ref{eq:xtoinf})} \\ \hline
       & number of terms $k$ & range of $x$ & number of terms $k$ & range of $x$ \\ \hline
[0.1, 0.2]    & $k=1$  & $x<10^{-16}$ &  & \\ \cline{1-3}
(0.2, 0.5]    &      & $x<10^{-8}$  & $k=10$ & $x> 10^{\frac{3}{1+\alpha}}$ \\ \cline{1-1}\cline{3-3}
(0.5, 0.99]   & \raisebox{.7em}[0pt][0pt]{$k=5$}  & $x<10^{-5}$  &  &  \\ \hline
(0.99,1.01]& \multicolumn{4}{c|}{formula (\ref{eq:alpha-1})} \\ \hline
(1.01,1.99999]& $k=10$ & $x<10^{-5}$  
& $k=10$ &$x> 10^{\frac{3}{1+\alpha}}$\\ \hline
(1.99999,2.0] & $k=85$ & $x\leq 7$   & \multicolumn{2}{c|}{
$\max((\ref{dense}),(\ref{eq:alpha-2}))$} \\ \hline
\end{tabular}
\label{tdense}
\end{center}
\end{table}

\section{Partial derivatives of the symmetric stable density with
  respect to the location and the scale parameters}
\label{sec:derivative-location-scale}

In this section we discuss the first and the second derivatives of the
stable density with respect to the location parameter $\mu$ and the
scale parameter $\sigma$. 

\subsection{The first derivatives w.r.t.\ location and scale}
\label{sec:first-derivative-location-scale}

Differentiating (\ref{dense}) we obtain
\begin{equation}
 f'(x;\alpha)=-\frac{1}{x}f(x;\alpha)+\frac{\alpha^2 \mbox{sign} (\alpha-1)}{\pi
 x^2 (\alpha-1)^2}
 \int_0^{\frac{\pi}{2}}
 g(\varphi;\alpha,x)(1-g(\varphi;\alpha,x))\exp(-g(\varphi;\alpha,x)) d\varphi
\label{eq:first-differential-dense}
\end{equation}
for $\alpha\neq 1$ and $x>0$. 
In view of (\ref{eq:f-mu}) and (\ref{eq:f-sigma})
we only have to evaluate (\ref{eq:first-differential-dense}).
The integrand $g(\cdot)(1-g(\cdot))\exp(-g(\cdot))$ 
has a positive local maximum  and a negative
local minimum.
These are attained at
$\varphi_2$: $g(x;\alpha,\varphi_2)=\frac{3-\sqrt{5}}{2}$ (local maximum) and
$\varphi_3$: $g(x;\alpha,\varphi_3)=\frac{3+\sqrt{5}}{2}$ (local minimum).
For $\alpha<1$ the ordering of these points is
$\varphi_3 \leq \varphi_1 \leq \varphi_2$
and the reverse
holds for $\alpha>1$. Accordingly we divide the interval of integration
and then the method of adaptive integration gives reliable values.
Without these divisions the numerical integration sometimes does not
converge or produces incorrect values. Note that
if $\alpha$ is very close to  $2$, the values of
$\varphi_i$, $i=1,2,3$, are  outside  of the integration
range. Therefore another treatment is needed for this case.

As in the case of the density itself, we use alternative representations of the first derivative
of the symmetric stable density around the boundary and following
representations are needed. 

\begin{enumerate}
\setlength{\itemsep}{2pt}
\item 
Expansions of the first derivative of the density.
\begin{equation}
f'(x;\alpha)=
\frac{1}{\pi\alpha}\sum^{\infty}_{k=1}\frac{\Gamma((2k+1)/\alpha)}{(2k-1)!}
(-1)^k x^{2k-1}.
\label{eq:expansion-first-differential-to-0} 
\end{equation}
As in (\ref{eq:xto0}), for $1 < \alpha \le 2$, this
series is convergent for every $x$.
For $\alpha < 1$, this series is justified as an asymptotic expansion as
$x\rightarrow 0$.  Similarly 
for the case $x\to \infty$  we have
\begin{equation}
f'(x;\alpha)= \frac{1}{\pi}\sum^{\infty}_{k=1}\frac{\Gamma(\alpha k+2)}{k!} 
 (-1)^k \sin(\frac{\pi \alpha k}{2})x^{-k\alpha-2}. 
\label{eq:expansion-first-differential-to-inf} 
\end{equation}
For $\alpha < 1$ this series converges for every $x\neq 0$ and 
for $\alpha>1$ this series is justified as an asymptotic expansion.

\item Taylor expansion of derivative of density around $\alpha=1$:
\begin{equation}
f'(x;\alpha) = f'(x;1) + f'_\alpha(x;1)(\alpha-1) 
+ \frac{1}{2}f'_{\alpha\alpha}(x;1)(\alpha-1)^2+ 
o((\alpha-1)^2),
\label{DTAEC}
\end{equation}
where
$$
f'(x;1)=-\frac{2x}{\pi(1+x^2)^2},
$$
\begin{eqnarray*}
f'_\alpha(x;1)
&=& 
\frac{1}{\pi}
\left\{\frac{-2x^3+6x}{(1+x^2)^3}\left(\frac{3}{2}-\gamma-\frac{1}{2}\log(1+x^2)\right)
+\frac{2-6x^2}{(1+x^2)^3}\arctan x\right\}, 
\label{fisrt term of TAEDC}
\end{eqnarray*}
and 
\begin{eqnarray*} \quad
f'_{\alpha\alpha}(x;1) &=& 
-\frac{2x(x^4-14x^2+9)}{\pi(1+x^2)^4}
\left\{\frac{\pi^2}{6}+\left(1-\gamma-\frac{1}{2}\log(1+x^2)\right)^2-1-\arctan^2x\right\} \\
& &-
\frac{8(3x^4-8x^2+1)}{\pi(1+x^2)^4}\arctan x
\left(\frac{3}{2}-\gamma-\frac{1}{2}\log(1+x^2)\right)
\nonumber   \\
& &-
\frac{2x(x^4-22x^2+17)}{\pi(1+x^2)^4}
\left(1-\gamma-\frac{1}{2}\log(1+x^2)\right)
\nonumber    \\
& &-
\frac{4(x^4-6x^2+1)}{(1+x^2)^4}\arctan x
+
\frac{8x(x^2-1)}{(1+x^2)^4} .
\end{eqnarray*}
\end{enumerate}

\begin{table}[htbp]
\caption[percentage]{Approximations to $f'(x;\alpha)$ at boundary cases} 
\label{tbl:first-differential-dense}
\vspace{-3mm}
\begin{center}
\begin{tabular}{|c|c|c|c|c|} \hline
$\alpha\backslash x$  & \multicolumn{2}{c|}{$x \to 0$, formula (\ref{eq:expansion-first-differential-to-0})  } 
& \multicolumn{2}{c|}{$x \to \infty$, formula (\ref{eq:expansion-first-differential-to-inf}) }  \\ \hline
   & number of terms $k$ & range of $x$ & number of terms $k$ & range of $x$ \\ \hline
[0.2, 0.25]  &      & $x<10^{-8}$  &  & \\ \cline{1-1}\cline{3-3}
(0.25, 0.3]  & $k=5$  & $x<10^{-6}$  & $k=10$ & $x> 10^{\frac{3}{1+\alpha}}$ \\ \cline{1-1}\cline{3-3}
(0.3, 0.99]  &      & $x<10^{-5}$  &  &  \\ \hline
(0.99,1.01]  & \multicolumn{4}{c|}{formula (\ref{DTAEC}) } \\ \hline
(1.01,1.99999]& \raisebox{-.7em}[0pt][0pt]{$k=10$} & \raisebox{-.7em}[0pt][0pt]{$x<10^{-3}$}  & $k=10$ &$x> 10^{\frac{3}{1+\alpha}}$\\ 
\cline{1-1}\cline{4-5}
(1.99999,2.0]  &     &   &
\multicolumn{2}{c|}{formula (\ref{eq:first-differential-dense})}\\ \hline
\end{tabular}
\label{tddense}
\end{center}
\end{table}

Table \ref{tddense} is a summary of approximations to the first
derivative $f'(x;\alpha)$.
The interpretation of Table \ref{tddense} is almost the same as Table
\ref{tdense}.
Note that for $0.99 < \alpha \le 1.01$ 
the approximation (\ref{DTAEC}) around Cauchy using
Taylor expansion to the order $O((\alpha-1)^2)$ 
is accurate enough.  
In the range $\alpha \in 
[0.1,0.2)$, approximations are not good and we do not consider
accurate calculations of derivatives. 
For $\alpha$ very close to 2 and $x$ large, we have the same problem
as in the case of the density.  However 
(\ref{eq:first-differential-dense}) gives reasonable values and we use
(\ref{eq:first-differential-dense}).
Note that if $\alpha$ is very close to
2, then we do not observer very large $x$  and the approximation to
$f'(x;\alpha)$ is less important  than the density itself. We will
discuss this problem again in Section 6 in connection with accurate
evaluation of the Fisher information for $\alpha$ very close to 2.

\begin{figure}[hbtp]
\begin{minipage}{.50\linewidth}
\includegraphics[width=\linewidth]{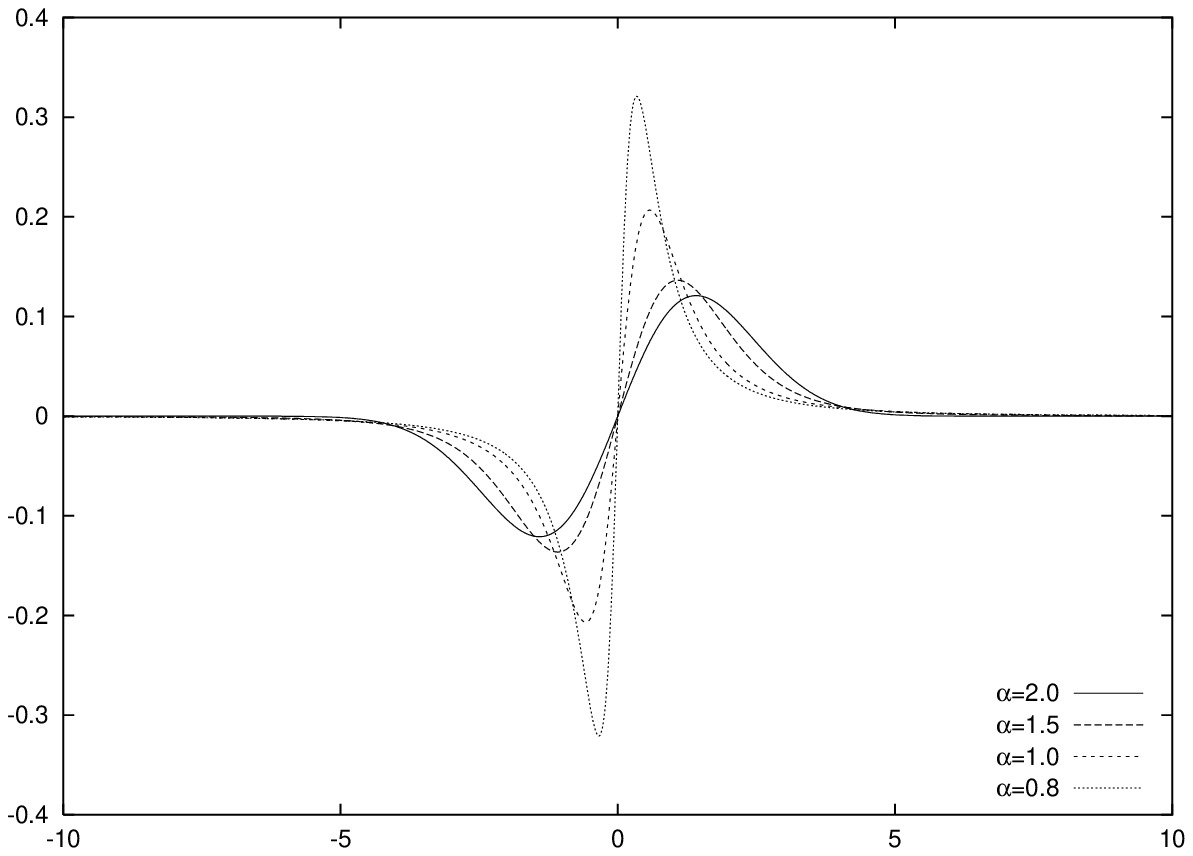}
\vspace{-5mm}
\caption{Derivative of density w.r.t.\  $\mu$}
\label{fig:dldense} 
\end{minipage}
\begin{minipage}{.50\linewidth}
\includegraphics[width=\linewidth]{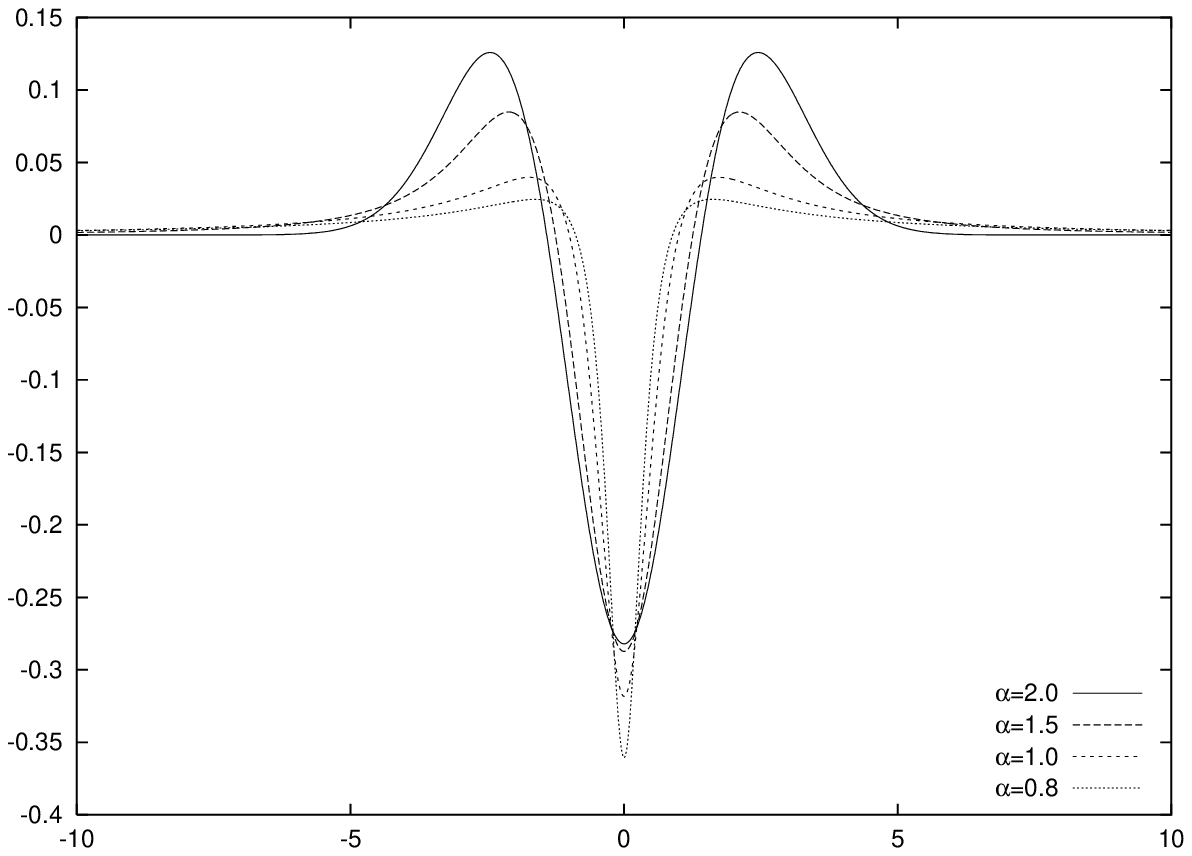}
\vspace{-5mm}
\caption{Derivative of density w.r.t.\ $\sigma$}
\label{fig:dsdense}
\end{minipage} 
\end{figure}

We present the graphs of the derivative concerning $\mu$ (Figure
\ref{fig:dldense}), and the derivative concerning $\sigma$ (Figure
\ref{fig:dsdense}).
Both derivatives are continuous in $\alpha$ and we do not see 
abrupt changes of the
forms for different values of $\alpha$.

\subsection{The second  derivatives w.r.t.\ location and scale}
\label{sec:second-derivative-location-scale}
We only need to evaluate 
$f''(x;\alpha)$ as in the case of the first derivatives.
Differentiating (\ref{eq:first-differential-dense}) we obtain
\begin{eqnarray}\label{eq:second-differential-dense}
 f''(x;\alpha) &=& \frac{1}{\alpha-1}\frac{1}{x^2}f(x;\alpha)
              +\frac{3-2\alpha}{\alpha-1}\frac{1}{x}f'(x;\alpha) \\
               && -\frac{\alpha^3 \mbox{sign} (\alpha-1)}{\pi
 (\alpha-1)^3 x^3}
 \int_0^{\frac{\pi}{2}}
 g^2(\varphi;\alpha,x)(2-g(\varphi;\alpha,x))\exp(-g(\varphi;\alpha,x)) d\varphi. \nonumber
\end{eqnarray}
Considering that $f(x;\alpha)$ and $f'(x;\alpha)$ can be calculated as
in the  previous sections, we only need 
to evaluate the integral which appear in the third
term of (\ref{eq:second-differential-dense}).
The integrand $g^2(\cdot)(2-g(\cdot))\exp(-g(\cdot))$ in the second
term on the right hand side of 
(\ref{eq:second-differential-dense}) has a zero at 
$\varphi_4$: $g(x;\alpha,\varphi_4)=2$, a local minimum 
at $\varphi_1$: $g(x;\alpha,\varphi_1)=1$
and a local maximum at $\varphi_5$: $g(x;\alpha,\varphi_5)=4$. For
$\alpha<1$ the order of these points is 
$\varphi_1 \leq \varphi_4 \leq \varphi_5$
and the reverse order holds for $\alpha>1$. 
Therefore we can divide the interval of integral according to these
values and get accurate results.
Alternative representations of the second derivative
of the symmetric stable density around the boundary
are as follows. 
\noindent
\begin{enumerate}
\setlength{\itemsep}{2pt}
\item 
Expansions of the second derivative of density.
\begin{equation}
f''(x;\alpha)=
\frac{1}{\pi\alpha}\sum^{\infty}_{k=1}\frac{\Gamma((2k+1)/\alpha)}{(2k-2)!}
(-1)^kx^{2k-2}, 
\label{eq:expansion-second-differential-to-0} 
\end{equation}
For $1<\alpha \leq 2$ this series is convergent for all $x$ 
and for $\alpha<1$ this is justified as an asymptotic expansion as $x \to 0$.
Similarly for $x \to \infty$ we have
\begin{equation}
f''(x;\alpha)=
 \frac{1}{\pi}\sum^{\infty}_{k=1}\frac{\Gamma(\alpha k+3)}{k!}
(-1)^{k-1} 
\sin(\frac{\pi\alpha k}{2}) x^{-k\alpha-3}, 
\label{eq:expansion-second-differential-to-inf} 
\end{equation}

For $\alpha<1$ this series is convergent for $x\neq 0$
and  for
$1 <\alpha \leq 2$ this series is justified as an asymptotic expansion.

\item Taylor expansion of the second derivative of density around $\alpha=1$:
\begin{equation}
f''(x;\alpha) = f''(x;1) + f''_\alpha(x;1)(\alpha-1) 
+ \frac{1}{2}f''_{\alpha\alpha}(x;1)(\alpha-1)^2+ 
o((\alpha-1)^2),
\label{eq:D2TAEC}
\end{equation}
where
$$
f''(x;1)=
\frac{6x^2-2}{\pi(1+x^2)^3},
$$
\begin{eqnarray*}
\label{eq:fdd_alpha}\quad 
f''_\alpha(x;1)
&=& 
\frac{2}{\pi}
\left\{\frac{x^4-6x^2+1}{(1+x^2)^4}\left(\frac{11}{2}-3\gamma-\frac{3}{2}\log(1+x^2)\right)
+\frac{12x(x^2-1)}{(1+x^2)^4}\arctan x\right\},
\end{eqnarray*}
and 
\begin{eqnarray*}
&&    \\
f''_{\alpha\alpha}(x;1) &=& 
\frac{6(x^6-25x^4+35x^3-3)}{\pi(1+x^2)^5}
\left\{\frac{\pi^2}{6}+\left(1-\gamma-\frac{1}{2}\log(1+x^2)\right)^2-1-\arctan^2x\right\} 
\nonumber \\
& &+
\frac{96x(x^4-5x^2+2)}{\pi(1+x^2)^5}\arctan x
\left(\frac{3}{2}-\gamma-\frac{1}{2}\log(1+x^2)\right)
\nonumber   \\
& &+
\frac{2(5x^6-155x^4+235x^2-21)}{\pi(1+x^2)^5}
\left(1-\gamma-\frac{1}{2}\log(1+x^2)\right)
\nonumber    \\
& &+
\frac{4x(11x^4-70x^2+31)}{\pi(1+x^2)^5}\arctan x
+ \frac{2(x^6-50x^4+85x^2-8)}{\pi(1+x^2)^5}.
\nonumber
\end{eqnarray*}
\end{enumerate}

\begin{table}[htbp]
\caption[percentage]{Approximations to $f''(x;\alpha)$ at boundary cases}
 \label{tbl:second-differential-dense}
\vspace{-3mm}
\begin{center}
\begin{tabular}{|c|c|c|c|c|} \hline
$\alpha \backslash  x$  & \multicolumn{2}{c|}{$x \to 0$,
formula (\ref{eq:expansion-second-differential-to-0})}
 & \multicolumn{2}{c|}{$x \to \infty$, 
formula (\ref{eq:expansion-second-differential-to-inf})}
 \\ \hline
  & number of terms $k$ & range of $x$ & number of terms $k$ & range of $x$ \\ \hline
[0.2, 0.25]  &      & $x<10^{-8}$  &  & \\ \cline{1-1}\cline{3-3}
(0.25, 0.3]  & \raisebox{-.7em}[0pt][0pt]{$k=5$}  & $x<10^{-6}$  
& \raisebox{-.7em}[0pt][0pt]{$k=10$} & 
\raisebox{-.7em}[0pt][0pt]{$x> 10^{\frac{3}{1+\alpha}}$}
 \\ \cline{1-1}\cline{3-3}
(0.3, 0.9]  &      & $x<10^{-5}$  &  &  \\ \cline{1-1}\cline{3-3}
(0.9, 0.99] &      & $x<10^{-3}$  &  &  \\ \hline
(0.99,1.01]  & \multicolumn{4}{c|}{formula (\ref{eq:D2TAEC}) } \\ \cline{1-1} \cline{2-3}
(1.01,1.02]  &     &    & \multicolumn{2}{c|}{} \\ \cline{1-1}  \cline{4-5}    
(1.02,1.999]& {$k=10$} 
& {$x<10^{-3}$}  & $k=10$ &$x> 10^{\frac{3}{1+\alpha}}$\\ \cline{1-1}\cline{4-5}
(1.999,2.0]  &  &   & \multicolumn{2}{c|}{formula 
(\ref{eq:second-differential-dense})} \\ \hline
\end{tabular}
\label{tddense2}
\end{center}
\end{table} 

Table \ref{tbl:second-differential-dense} summarizes uses of various
expansions. Here
formula (\ref{eq:D2TAEC}) means use of Taylor expansion of the second 
derivative of 
density around $\alpha=1$ and 
(\ref{eq:second-differential-dense})
means direct use of integral expression of the second derivative.
Note for $\alpha \to 2$ and $x \to \infty$ the same problem arise as in
the density. 
However we do not use the second
derivative of density for $\alpha \in (1.999,2.0]$ in this paper.

\section{Partial derivatives with respect to 
 the characteristic exponent} %
\label{sec:derivative-alpha}

In this section we discuss the first and the second derivatives of the
stable density with
respect to the characteristic exponent $\alpha$.

\subsection{The first derivative w.r.t. $\alpha$}
\label{sec:first-derivative-alpha}

We only need to investigate the standard ($\mu=0$ and $\sigma=1$) case.
Although the representations of $f_\alpha(x;\alpha)$ sometimes
become complicated, careful analysis of various cases  gives accurate and efficient
calculations. 
The direct differentiation of (\ref{dense}) yields the following
representation.
\begin{equation}
  f_\alpha(x;\alpha)=\frac{1}{\alpha(1-\alpha)}f(x;\alpha)
 +\frac{\alpha}{\pi |\alpha-1| x}
 \int_0^{\frac{\pi}{2}}
 g_\alpha(\varphi;\alpha,x)(1-g(\varphi;\alpha,x))\exp(-g(\varphi;\alpha,x))
 d\varphi,
\label{eq:first-differential-alpha}
\end{equation}
 where
\begin{eqnarray}
  g_\alpha(x;\alpha)
 &=&
-\left(\frac{x\cos \varphi}{\sin \alpha
 \varphi}\right)^{\frac{\alpha}{\alpha-1}}
   \frac{\cos(\alpha-1)}{\cos \varphi} 
\\ && \qquad\qquad \times
\left\{ \frac{1}{(\alpha-1)^2}\log\left|\frac{x\cos \varphi}{\sin \alpha \varphi} \right|+
       \frac{\alpha \varphi}{\alpha-1}\frac{1}{\tan \alpha\varphi}+ \varphi\tan(\alpha-1)\varphi\right\} \nonumber \label{diffag} \\
 &=& 
-g(x;\alpha)\left\{ \frac{1}{(\alpha-1)^2}\log\left|\frac{x\cos
  \varphi}{\sin \alpha \varphi} \right| 
  + \frac{\alpha \varphi}{\alpha-1}\frac{1}{\tan \alpha\varphi}
  + \varphi\tan(\alpha-1)\varphi\right\}
\nonumber \\
 &=& 
-g(x;\alpha) ( h_1(\varphi) + h_2(\varphi) + h_3(\varphi)).
\nonumber 
\end{eqnarray}
Here for convenience $h_i(\varphi)$, $i=1,2,3$, are defined as
\begin{equation}
h_1(\varphi)= \frac{1}{(\alpha-1)^2}\log\left|\frac{x\cos \varphi}{\sin \alpha \varphi} \right|,\quad h_2(\varphi)= \frac{\alpha \varphi}{\alpha-1}\frac{1}{\tan \alpha\varphi}, \quad 
h_3(\varphi)= \varphi\tan(\alpha-1)\varphi.
\end{equation}

We consider the integrand
$g_\alpha(\cdot)(1-g(\cdot))\exp(-g(\cdot))$ in (\ref{eq:first-differential-alpha})
separately, dividing it 
into three parts $g(\cdot)h_i(\cdot)(1-g(\cdot))\exp(-g(\cdot))$.

For each  $g(\cdot)h_i(\cdot)(1-g(\cdot))\exp(-g(\cdot))$, $i=1,2,3$, 
we divide the interval of integration appropriately based on zero
points $\varphi_i, \ i=1,4,5$.
The following properties are obtained by careful evaluation of the signs
of derivatives of $h_i(\cdot)$. 
\begin{description}
\item{$h_1(\varphi)$:} strictly decreasing, zero at
           $\varphi_4$: $\frac{x\cos \varphi}{\sin \alpha \varphi}=1$,
           $\varphi_1<\varphi_4$ for $\alpha <1$ and $\varphi_4<\varphi_1$
           for $\alpha>1$.
\item{$h_2(\varphi)$:} negative for $\alpha<1$. The sign
  changes from $+$ to $-$ at $\varphi_5= \frac{\pi}{2\alpha}$ for
           $\alpha>1$.
\item{$h_3(\varphi)$:} nonpositive for $\alpha<1$. nonnegative for $\alpha>1$.
\end{description}
As in the density, we use alternative representations of the first derivative
of the symmetric stable density concerning $\alpha$ around the boundary and following
representations are needed. 

\begin{enumerate}
\setlength{\itemsep}{2pt}
\item 
Expansions of the first derivative of the density w.r.t $\alpha$.
\begin{equation}
\label{eq:expansion-alpha1-zero}
f_\alpha(x;\alpha)=
-\frac{1}{\pi\alpha^2}\sum^{\infty}_{k=0}\frac{\Gamma'((2k+1)/\alpha+1)}{(2k)!}
(-x^2)^k.
\end{equation}
For $1 < \alpha \le 2$ this
series is convergent for every $x$ and
for  $\alpha < 1$ it is justified as an asymptotic expansion as
$x\rightarrow 0$.
Similarly for $x \to \infty$ we have
\begin{eqnarray}
\label{eq:expansion-alpha1-inf}
f_\alpha(x;\alpha)&=&
\frac{1}{\pi}\sum^{\infty}_{k=1}\frac{\Gamma'(\alpha k+1)}{(k-1)!}
(-1)^{k-1}\sin\left(\frac{\pi\alpha k}{2}\right)x^{-k\alpha-1}
\\ && 
+\frac{1}{\pi}\sum^{\infty}_{k=1}\frac{\Gamma(\alpha k+1)}{(k-1)!}
(-1)^{k-1}
\left[\frac{\pi}{2}\cos\left(\frac{\pi\alpha k}{2}\right)
   - \log x\; \sin\left(\frac{\pi\alpha k}{2}\right)
     \right]x^{-k\alpha-1}.  
\nonumber 
\end{eqnarray}
For $\alpha < 1$ this series converges for every $x\neq 0$ and 
for $\alpha>1$ this series is justified as an asymptotic expansion.

\item 
Taylor expansion of the first derivative of $\alpha$ around $\alpha=1$:
\begin{equation}
f_{\alpha}(x;\alpha) \label{DATEAC}
= f_\alpha(x;1)+ f_{\alpha\alpha}(x;1)(\alpha-1) 
+ o(|\alpha-1|),
\end{equation} 
where $f_\alpha(x;1)$ and $f_{\alpha\alpha}(x;1)$ are
given by (\ref{eq:alpha-1-1}) and (\ref{eq:alpha-1-2}).
\end{enumerate}

\begin{figure}
\begin{minipage}{.50\linewidth}
\includegraphics[width=\linewidth]{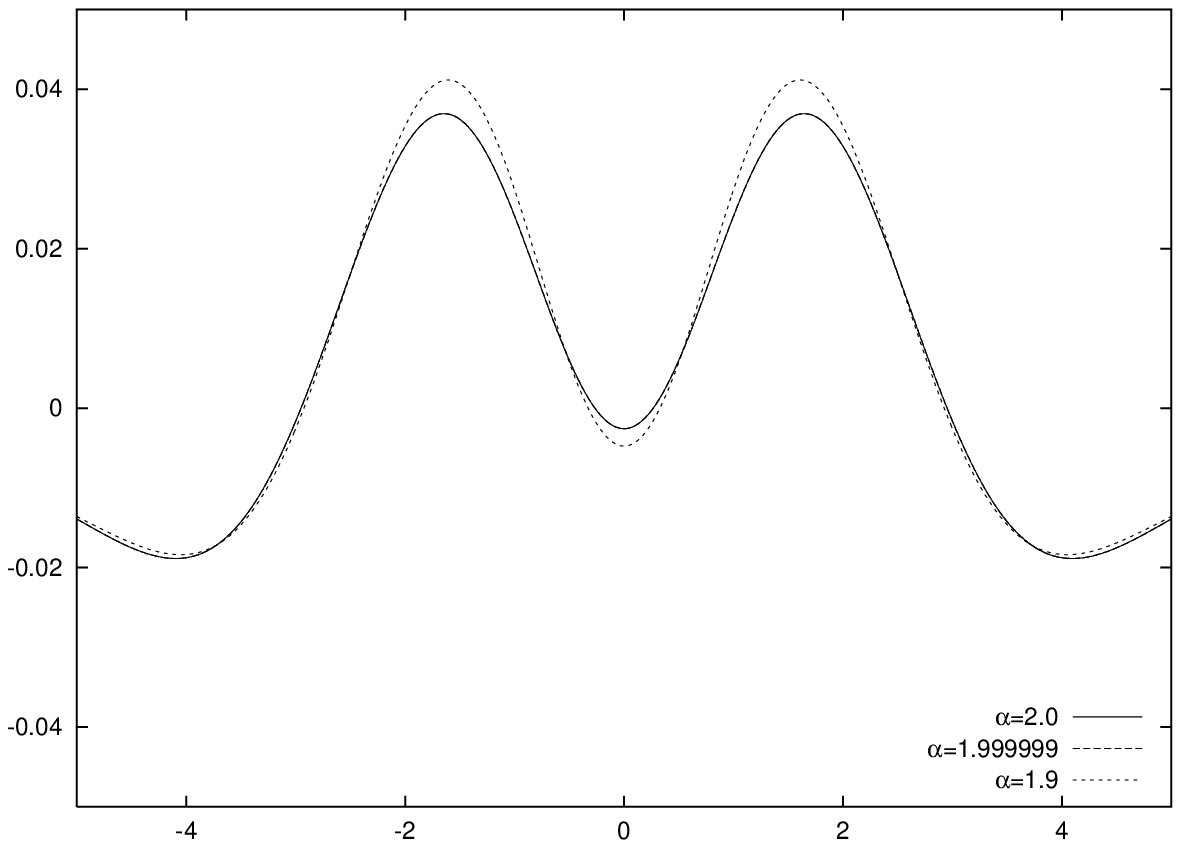}
\vspace{-5mm}
\caption{$f_\alpha(x;\alpha)$ (accurate)}
\label{fig:dadense3} 
\end{minipage}
\begin{minipage}{.50\linewidth}
\includegraphics[width=\linewidth]{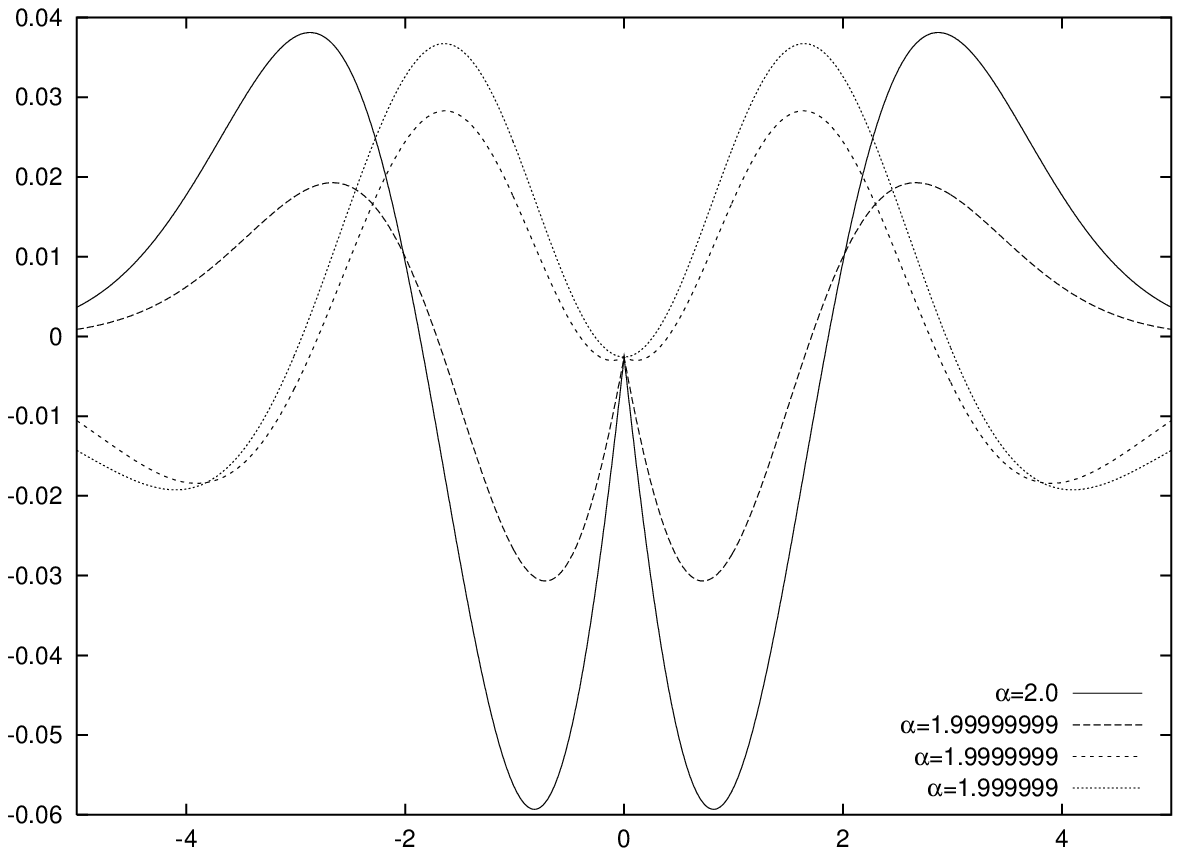}
\vspace{-5mm} 
\caption{$f_\alpha(x;\alpha)$ using 
(\ref{eq:first-differential-alpha}) (not accurate)}
\label{fig:dadense4}
\end{minipage} 
\end{figure}

In order to first confirm the derivative at $\alpha=2$ we utilize
equation 3.21 on page 12 of Oberhettinger (1990).
\begin{equation*}
\int_0^\infty  u^{\nu}e^{-a u^2}\cos(u y)du = 
\frac{1}{2}\left(\frac{\pi}{2}\right)^{\frac{1}{2}}
(2a)^{-\frac{1}{2}-\frac{1}{2}\nu}
\sec\left(\frac{\pi\nu}{2}\right)\exp\left(-\frac{y^2}{8a}\right)
\left(D_\nu(z)+D_\nu(-z)\right),
\end{equation*}
where
$z=(2a)^{-\frac{1}{2}}y$ %
and $D_\nu(z)$ is the parabolic cylinder functions (see p.255 of Oberhettinger (1990)). 
Differentiating this with respect to $\nu$, setting $\nu =2$  and
adjusting some signs and constants we obtain $f_{\alpha}(x;2)$. 
These values coincide with
that of expansions in 
(\ref{eq:expansion-alpha1-zero})
at $\alpha=2$.
Note the integral (\ref{eq:first-differential-alpha}) is not accurate
when $\alpha$ is very close to $2$ especially for 
$\alpha >1.9999999$ . This inaccuracy of  the integral
(\ref{eq:first-differential-alpha}) is illustrated in  Figure
\ref{fig:dadense4}, 
in comparison 
to the accurate  values shown in Figure \ref{fig:dadense3}.

Table \ref{tbl:ttdadense} is a summary of approximations to the first 
derivative
of the density function with respect to $\alpha$. 

The formula (\ref{DATEAC}) means Taylor expansion of the first derivative of
$\alpha$ around $\alpha=1$. For $\alpha \in (1.9999,2.0]$ we use only
the expansion (\ref{eq:expansion-alpha1-zero}) and the 
asymptotic expansion (\ref{eq:expansion-alpha1-inf}).
   
\begin{table}[htbp]
\caption[percentage]{Approximations to $f_\alpha(x;\alpha)$ at
  boundary cases}
\label{tbl:ttdadense}
\vspace{-3mm}
\begin{center}
\begin{tabular}{|c|c|c|c|c|} \hline
$\alpha \backslash x$  & \multicolumn{2}{c|}{$x \to 0$,
formula(\ref{eq:expansion-alpha1-zero})
} & \multicolumn{2}{c|}{$x \to \infty$,
formula(\ref{eq:expansion-alpha1-inf})
}  \\ \hline
 & number of terms $k$ & range of $x$ 
& number of terms $k$ & range of $x$ \\ \hline
[0.1, 0.2]   & $k=1$  & $x<10^{-16}$ &  & \\ \cline{1-3}
(0.2, 0.3]   &    & $x<10^{-7}$  
& \raisebox{-.7em}[0pt][0pt]{$k=10$}  
& \raisebox{-.7em}[0pt][0pt]{$x> 10^{\frac{3}{1+\alpha}}$} \\ \cline{1-1}\cline{3-3}
(0.3, 0.5]   & $k=5$  & $x<10^{-5}$  &  & \\ \cline{1-1}\cline{3-3}
(0.5, 0.99]  &  & $x<10^{-5}$  &  & \\ \hline
(0.99,1.01]& \multicolumn{4}{c|}{formula (\ref{DATEAC}) } \\ \hline
(1.01,1.9999]& $k=10$ 
& $x<10^{-5}$  & $k=10$ &$x> 10^{\frac{3}{1+\alpha}}$\\ \hline
(1.9999,2.0]  & $k=85$ & $x \le 8.0$ & $k=20$ & $x>8.0$ \\ \hline
\end{tabular}
\label{tdadense} 
\end{center}
\end{table}

\begin{figure}
\begin{minipage}{.50\linewidth}
\includegraphics[width=\linewidth]{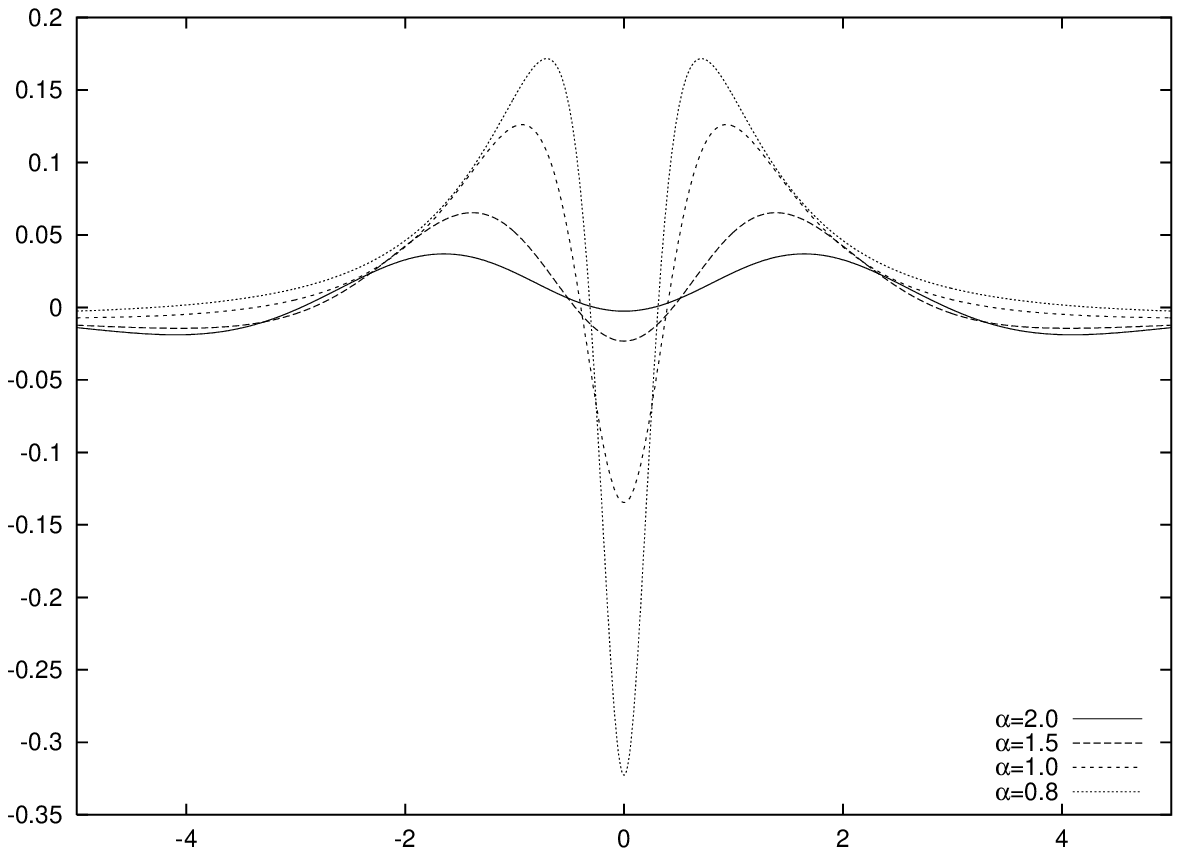}
\vspace{-5mm}
\caption{First derivative $f_\alpha(x;\alpha)$ w.r.t.\ to $\alpha$}
\label{fig:dadense1} 
\end{minipage}
\begin{minipage}{.50\linewidth}
\includegraphics[width=\linewidth]{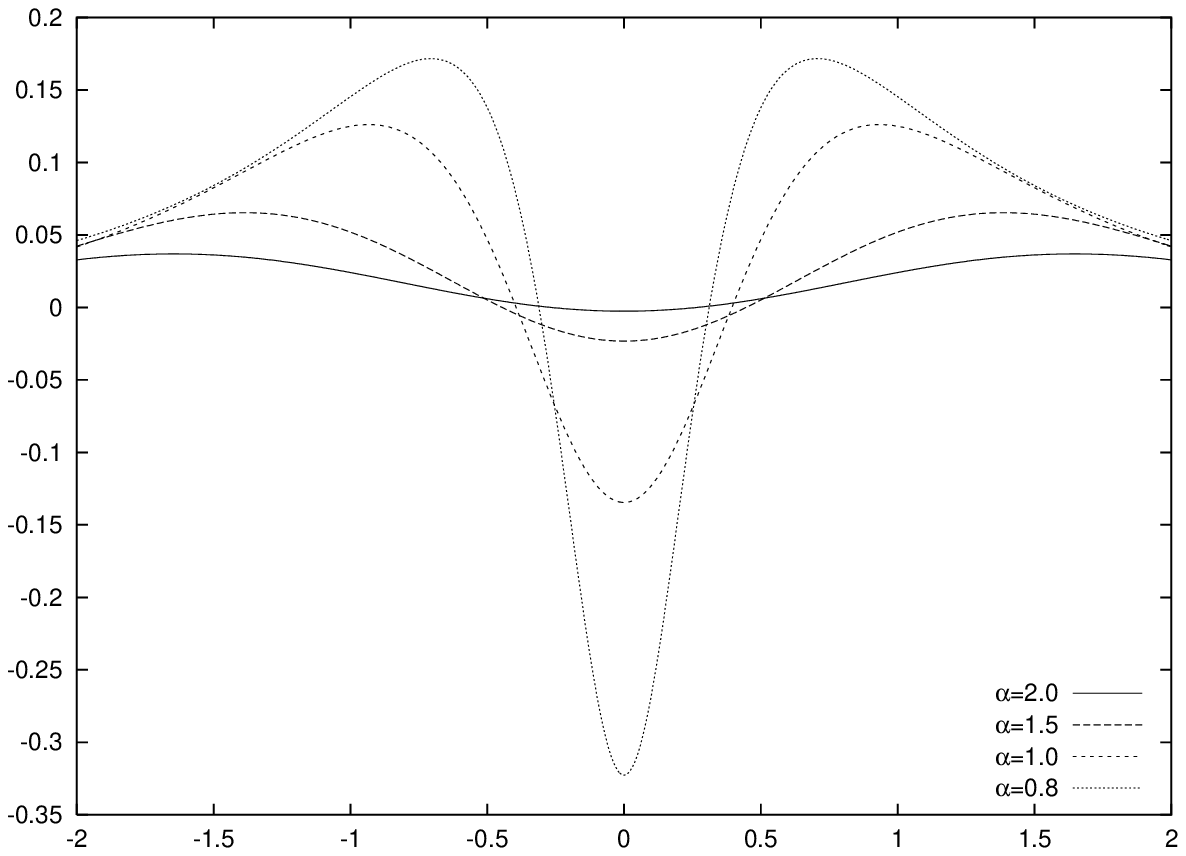}
\vspace{-5mm}
\caption{$f_\alpha(x;\alpha)$ (more detail around 0)}
\label{fig:dadense2} 
\end{minipage}
\end{figure}

\begin{figure}
\begin{minipage}{.50\linewidth}
\includegraphics[width=\linewidth]{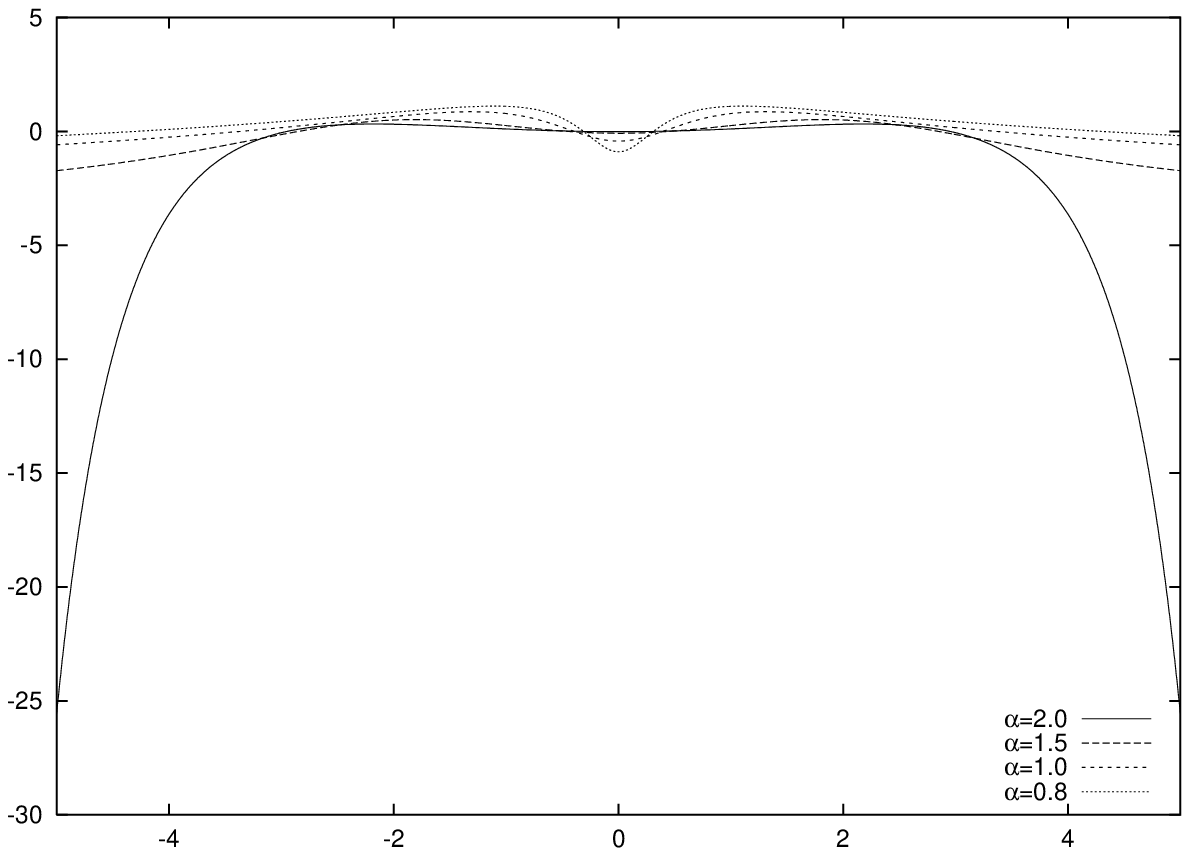}
\vspace{-5mm}
\caption{Score function $f_\alpha(x;\alpha)/f(x;\alpha)$}
\label{fig:alphascore1}
\end{minipage} 
\begin{minipage}{.50\linewidth}
\includegraphics[width=\linewidth]{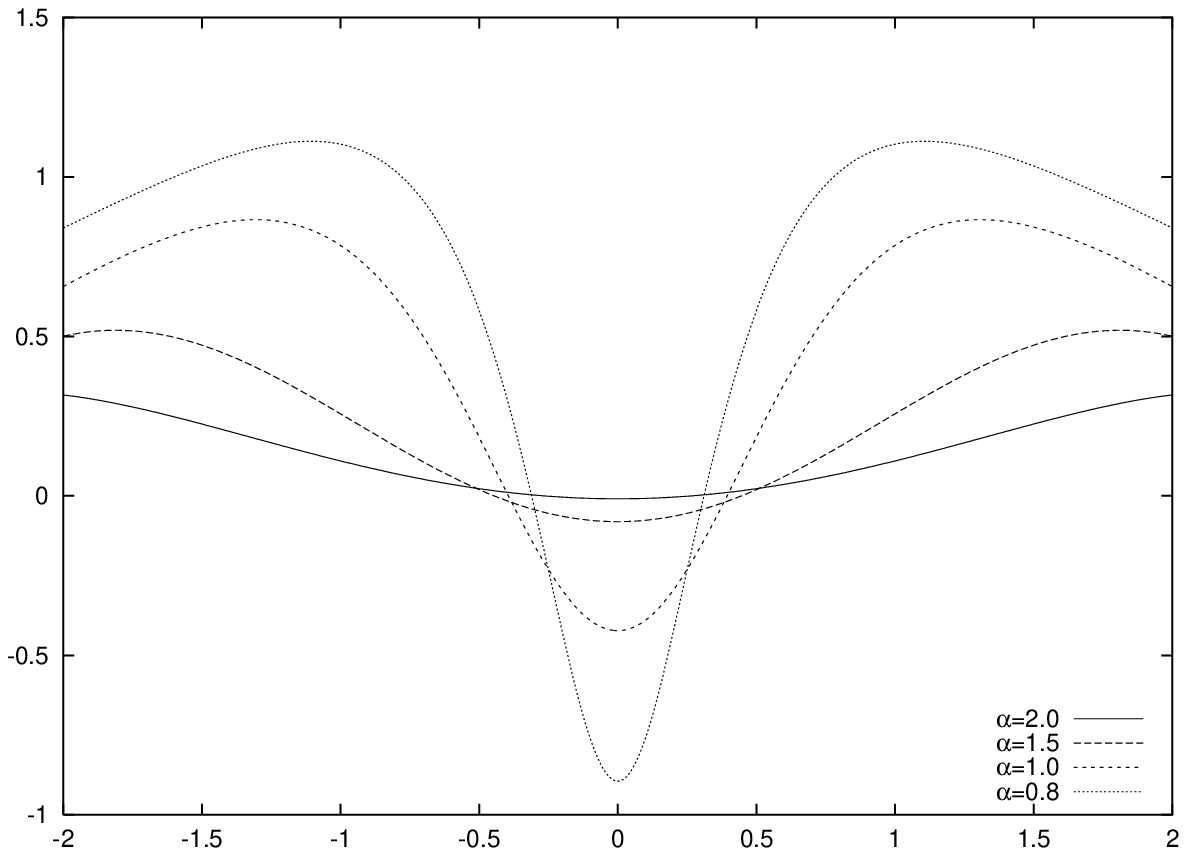}
\vspace{-5mm}
\caption{$f_\alpha(x;\alpha)/f(x;\alpha)$ (more detail around 0)}
\label{fig:alphascore2}
\end{minipage} 
\end{figure}

Here we present graphs of the first derivative w.r.t.\ $\alpha$ in Figure
\ref{fig:dadense1}, Figure \ref{fig:dadense2}, 
Figure \ref{fig:dadense3}, Figure \ref{fig:dadense4},
Figure \ref{fig:alphascore1} and Figure \ref{fig:alphascore2}. 
The first derivative $f_\alpha(x;\alpha)$ with respect to $\alpha$ is 
shown in Figure \ref{fig:dadense1} and \ref{fig:dadense2}
and  the score function $f_\alpha(x;\alpha)/f(x;\alpha)$ is shown in 
Figure \ref{fig:alphascore1} and \ref{fig:alphascore2}.
Note that graphs of $f_\alpha(x;\alpha)$ where $\alpha=1.999999$ and
$\alpha=2.0$ in Figure \ref{fig:dadense3} are almost the same
corresponding to the fact $f_\alpha(x;\alpha)$ is continuous in 
$\alpha \in (0,2]$. %

\subsection{The second derivative w.r.t.\ $\alpha$}
\label{sec:second-derivative-alpha}
Differentiating (\ref{eq:first-differential-alpha}) we obtain
\begin{eqnarray}
&&  f_{\alpha\alpha}(x;\alpha) = \label{eq:second-differential-alpha}
  \frac{2}{\alpha(1-\alpha)^{2}}f(x;\alpha)   \\
&& \qquad  +\frac{2\sgn(\alpha-1)}{\pi (\alpha-1)^{2}x}
   \int_0^{\frac{\pi}{2}}
    g(\varphi)\left\{(1+\alpha)h_1(\varphi)+2h_2(\varphi)+h_3(\varphi)\right\}(1-g(\varphi))
    \exp(-g(\varphi))d\varphi  \nonumber  \\
&&\qquad + \frac{\alpha}{\pi|\alpha-1|x}\int^{\frac{\pi}{2}}_0 
     g(\varphi)\left(g^2(\varphi)-3g(\varphi)+1\right)
     \left\{h_1(\varphi)+h_2(\varphi)+h_3(\varphi)\right\}^{2}
     \exp(-g(\varphi)) d\varphi \nonumber    \\
&&\qquad + \frac{\alpha}{\pi|\alpha-1|x}\int^{\frac{\pi}{2}}_0 
     g(\varphi)
    \left(\frac{\alpha}{\alpha-1}\frac{\varphi^2}{\sin^2\alpha\varphi}
           -\frac{\varphi^2}{\cos^2(\alpha-1)\varphi}\right)(1-g(\varphi))\exp(-g(\varphi)) d\varphi. \nonumber 
\end{eqnarray}
Concerning the second integral
we utilize the arguments of Section 
\ref{sec:first-derivative-location-scale} for efficient
integration. The third
integral is calculated by dividing the interval of integration
using zeros of integrand $\varphi_i$ for $i=1,2,3$.
As for the fourth integral, we integrate separately the two terms in the %
integrand. 
An alternative representations of the second derivative
of the symmetric stable density concerning $\alpha$ around the boundary
are as follows.

\begin{enumerate}
\setlength{\itemsep}{2pt}
\item 
Expansions of the second derivative of density w.r.t. $\alpha$.
\begin{eqnarray}
\label{eq:expansion-alpha2-zero}
f_{\alpha\alpha}(x;\alpha)& =&
\frac{2}{\pi\alpha^3}\sum^{\infty}_{k=0}\frac{\Gamma'((2k+1)/\alpha+1)}{(2k)!}
(-x^2)^{k} 
\\ 
&&
+\frac{1}{\pi\alpha^4}\sum^{\infty}_{k=0}\frac{\Gamma''((2k+1)/\alpha+1)}{(2k)!}
(2k+1)(-x^2)^{k}.
\nonumber
\end{eqnarray}
For $1<\alpha \leq 2$ this series is convergent for all $x$ 
and for $\alpha<1$ this is justified as an asymptotic expansion as $x \to 0$.
Similarly for $x \to \infty$ we have
\begin{eqnarray}
\label{eq:expansion-alpha2-inf}
f_{\alpha\alpha}(x;\alpha)& =&
\frac{1}{\pi}\sum^{\infty}_{k=1}\frac{\Gamma''(\alpha k+1)}{(k-1)!}
k  
(-1)^{k-1} \sin\left(\frac{\pi\alpha k}{2}\right) x^{-k\alpha-1}  
\\   &&
-\frac{2}{\pi}\sum^{\infty}_{k=1}\frac{\Gamma'(\alpha k+1)}{(k-1)!}
k
(-1)^{k} \left[\frac{\pi}{2}\cos\left(\frac{\pi\alpha k}{2}\right)
  - \log  x\; \sin\left(\frac{\pi\alpha k}{2}\right) \right]
x^{-k\alpha-1}
\nonumber  \\  &&
+\frac{1}{\pi}\sum^{\infty}_{k=1}\frac{\Gamma(\alpha k+1)}{(k-1)!}
k (-1)^{k-1}\left[\Big( \log^2(x)- \frac{\pi^2}{4}\Big) 
  \sin\left(\frac{\pi\alpha k}{2}\right)
\right.  
\nonumber \\ &&  \qquad
- \left. 
\pi\log x\; \cos\left(\frac{\pi\alpha k}{2}\right) \right] x^{-k\alpha-1}.
\nonumber
\end{eqnarray}
For $\alpha<1$ this series is convergent for $x\neq 0$
and  for
$1 <\alpha \leq 2$ this series is justified as an asymptotic expansion.

\item  Taylor expansion of the second derivative of $\alpha$ around $\alpha=1$:
\begin{equation}
f_{\alpha\alpha}(x;\alpha) \label{eq:DA2TEAC}
= f_{\alpha\alpha}(x;1)+ f_{\alpha\alpha\alpha}(x;1)(\alpha-1) 
+ o(|\alpha-1|),
\end{equation} 
where $f_{\alpha\alpha\alpha}(x;1)$ is given in Appendix \ref{sec:alpha-3}.
\end{enumerate}

Table \ref{tbl:second-differential-alpha} is a summary of
approximations of the second derivative of density function concerning
$\alpha$. Formula (\ref{eq:DA2TEAC}) means use of Taylor expansion of
the second derivative of $\alpha$ around $\alpha=1$.  For the sake of
convenience we use the integral representation
(\ref{eq:second-differential-alpha}) when $\alpha \in (1.999,2.0)$ and
$x \to \infty$.  This is somewhat problematic as in the case of the
first derivative w.r.t.\ $\alpha$ and further investigation is needed
for $\alpha$ near $2$.  In the calculations of the present paper, this
problem appears only in the observed fisher information for
$\alpha=1.5$. However the estimated values $\hat{\alpha}$
very near $\alpha=2$ seldom occur and there seems to be no influence of
this problem in the values presented in Table \ref{tbl:observed-info-alpha}.

\begin{table}[htbp]
\caption[percentage]{Approximations to $f_{\alpha\alpha}(x;\alpha)$ at
  boundary cases}
\label{tbl:second-differential-alpha}
\vspace{-3mm}
\begin{center}
\begin{tabular}{|c|c|c|c|c|} \hline
$\alpha \backslash x$  
& \multicolumn{2}{c|}{$x \to 0$,
formula (\ref{eq:expansion-alpha2-zero})
} & \multicolumn{2}{c|}{$x \to \infty$,
formula (\ref{eq:expansion-alpha2-inf})
}  \\ \hline
    & number of terms $k$ & range of $x$ & number of terms $k$ & range of $x$ \\ \hline
[0.2, 0.3]   &    & $x<2.0\times10^{-7}$  
& &  \\ \cline{1-1}\cline{3-3}
(0.3, 0.5]   & $k=5$  & $x<10^{-5}$  
& {$k=10$}  
& {$x> 10^{\frac{3}{1+\alpha}}$}
\\ \cline{1-1}\cline{3-3}
(0.5, 0.99]  &  & $x<10^{-3}$  &  & \\ \hline
(0.99,1.03]& \multicolumn{4}{c|}{formula (\ref{eq:DA2TEAC}) } \\ \hline
(1.01,1.999]& \raisebox{-.7em}[0pt][0pt]{$k=10$} 
& \raisebox{-.7em}[0pt][0pt]{$x<10^{-3}$}  & $k=10$ &$x> 10^{\frac{3}{1+\alpha}}$\\ \cline{1-1}\cline{4-5}
(1.999,2.0)  & & & \multicolumn{2}{c|}{formula (\ref{eq:second-differential-alpha})}  \\ \hline
\end{tabular}
\label{tdadense2} 
\end{center}
\end{table}

\section{Fisher informations of symmetric stable distributions}
\label{sec:fisher-information}

Using the results of the previous section we can accurately
evaluate the Fisher information matrices  of the symmetric stable
distributions.
To the
authors' knowledge there are only two somewhat incomplete results concerning
the Fisher Informations of stable distributions. One is Nolan (2001) and
the other is DuMouchel (1975).  The former was based on 
numerical differentiations of densities and did not give the
informations for $\alpha < 0.5$. The latter had the problem of
truncation of the integral 
and did not give the informations for $\alpha < 0.8$. Considering that
the tails and modes of distributions are important for calculating
informations, our approach is desirable.  We give the informations
$\alpha \in [0.2,2.0]$ of the symmetric distributions in Table
\ref{tbl:information}.  In Table \ref{tbl:information} `$\ast$' means
that the quantity is not defined.
The values  in Table \ref{tbl:information}
largely coincide with those of DuMouchel (1975). 

The components of Fisher information matrix $I$ are defined as follows.
\begin{equation}
\label{eq:information}
   I_{ij}=\int_{-\infty}^{\infty} \frac{\partial f}{\partial \theta_i}
\frac{\partial f}{\partial \theta_j} \frac{1}{f} dx.
\end{equation}
where $(\theta_1, \theta_2, \theta_3)=(\mu, \sigma,\alpha)$.  By
symmetry it can be easily shown that $I_{12}=I_{21}=0$.  For
notational clarity we write $I_{11}=I_{\mu\mu}$,
$I_{22}=I_{\sigma\sigma}$, $I_{23}=I_{\sigma\alpha}$ and
$I_{33}=I_{\alpha\alpha}$.  For numerical calculations of
(\ref{eq:information})
we use adaptive integration of infinite intervals with
singularities (QUAGIU) in GNU Scientific Library (2003). In can be
shown that the information matrix at the Cauchy ($\alpha=1.0$)
is analytically given as
\begin{equation}
I_{\mu\mu}=I_{\sigma\sigma}=1/2\sigma^2,
\quad I_{\sigma\alpha}=
\frac{1}{2\sigma}(1-\gamma-\log2),\quad
 I_{\alpha\alpha}=\frac{1}{2}
\left\{\frac{\pi^2}{6}+{(\gamma+\log2-1)}^2\right\}.
\end{equation}
These values are consistent with numerical calculations in Table \ref{tbl:information}.

\begin{table}[htbp]
\caption[percentage]{Fisher information matrix of symmetric stable distributions}
\begin{center}
\begin{tabular}{ccccc}    \hline
$\alpha \setminus I_{ij} $   & $I_{\mu\mu}$ & $I_{\sigma\sigma}$ & $I_{\alpha\alpha}$ & $I_{\sigma\alpha}$ \\ \hline
2.0     &0.5 & 2.0 & $ \infty $     & $\ast $ \\      
1.999   & 0.4995 & 1.9904 & 29.461  & $-$0.8685 \\
1.99    & 0.4960 & 1.9321 & 4.6197  & $-$0.6682 \\
1.95    & 0.4842 & 1.7631 & 1.4108  & $-$0.4821 \\
1.9     & 0.4727 & 1.6127 & 0.8846  & $-$0.3963 \\
1.8     & 0.4552 & 1.3898 & 0.5937  & $-$0.3138 \\
1.7     & 0.4424 & 1.2189 & 0.5028  & $-$0.2692 \\
1.6     & 0.4334 & 1.0775 & 0.4726  & $-$0.2396 \\
1.5     & 0.4281 & 0.9556 & 0.4737  & $-$0.2174 \\ 
1.4     & 0.4270 & 0.8475 & 0.4973  & $-$0.1992 \\
1.3     & 0.4310 & 0.7498 & 0.5424  & $-$0.1832 \\
1.2     & 0.4419 & 0.6603 & 0.6119  & $-$0.1679 \\
1.1     & 0.4630 & 0.5774 & 0.7132  & $-$0.1523 \\
1.05    & 0.4790 & 0.5381 & 0.7794  & $-$0.1440 \\
1.01    & 0.4953 & 0.5075 & 0.8413  & $-$0.1369 \\
1.0     & 0.5    & 0.5    & 0.8590  & $-$0.1352 \\
0.99    & 0.5049 & 0.4925 & 0.8763  & $-$0.1332 \\
0.95    & 0.5276 & 0.4631 & 0.9552  & $-$0.1257 \\
0.9     & 0.5641 & 0.4272 & 1.0721  & $-$0.1154 \\
0.8     & 0.6800 & 0.3586 & 1.3928  & $-$0.0913 \\
0.7     & 0.9094 & 0.2937 & 1.8974  & $-$0.0611 \\
0.6     & 1.4446 & 0.2325 & 2.7414  & $-$0.0220 \\
0.5     & 3.1167 & 0.1753 & 4.2748  &  0.0295 \\
0.4     & 12.256 & 0.1226 & 7.3994  &  0.0979 \\
0.3     & 188.09 & 0.0756 & 14.925  &  0.1869 \\
0.2     & 149359.4 & 0.0367 & 38.729 & 0.2938 \\
\end{tabular}
\label{tbl:information}
\end{center}
\end{table}

In order to check our computations, 
we have done small simulation study of maximum
likelihood estimation of $\alpha$ (for fixed $\mu=0$ and $\sigma=1$).
We found that simulated information coincides well with the  exact
information. Table \ref{tbl:simlation} is the results of $1000$
iterations of maximum likelihood estimation
of $\alpha$ with the sample size of $n=50$. We denote the maximum
likelihood estimator by $\hat{\alpha}$.
$\bar{\alpha}$ is the average of $1000$ iterations
and  $\hat{\sigma}^2_\alpha$ is the variance of 
$\sqrt{50} \times (\hat{\alpha}-\alpha)$. Though the convergence of  
$1/\hat{\sigma}^2_{\alpha}$ to $I_{\alpha\alpha}$ seems somewhat slow for $\alpha<1$, 
the values of 
$1/\hat{\sigma}^2_{\alpha}$ are consistent with the exact Fisher information
$I_{\alpha\alpha}$. 

\begin{table}[htbp]
\caption[percentage]{MLE\ ($n=50$, 1000 iteration)}
\begin{center}
\begin{tabular}{cccc}    \hline
$\alpha$ & $\bar{\alpha}$ & $1/\hat{\sigma}_{\alpha}$ &  $I_{\alpha\alpha}$   \\ \hline
1.99     & 1.974  &  4.4877   &  4.6197    \\
1.8      & 1.811  &  0.6396   &  0.5937    \\
1.7      & 1.714  &  0.5102   &  0.5028    \\
1.5      & 1.531  &  0.4532   &  0.4737    \\
1.3      & 1.320  &  0.5114   &  0.5424    \\
1.0      & 1.027  &  0.7676   &  0.8590    \\
0.8      & 0.819  &  1.1944   &  1.3928    \\
0.5      & 0.509  &  3.8803   &  4.2748
\end{tabular}
\label{tbl:simlation}
\end{center}
\end{table}

Furthermore we simulated the observed Fisher information (\ref{eq:observed-FI})
of $\alpha$ %
in Table \ref{tbl:observed-info-alpha}.
Here $\hat{I}_{\alpha\alpha}(2)$
corresponds to the right hand side of (\ref{eq:observed-FI}) and 
$$
\hat{I}_{\alpha\alpha}(1)= \frac{1}{n}\sum_{i=1}^n \left(\frac{ 
  f_\alpha(x_i ; \hat\alpha)}{f(x_i ; \hat\alpha)}\right)^2
$$
involves the first derivative only.
Variance of $\hat{I}_{\alpha\alpha}(1)$ and 
$\hat{I}_{\alpha\alpha}(2)$ are also shown in the parentheses
in Table
\ref{tbl:observed-info-alpha}. $\hat{I}_{\alpha\alpha}(2)$ has smaller 
variance than $\hat{I}_{\alpha\alpha}(1)$ but some positive  bias is
observed in $\hat{I}_{\alpha\alpha}(2)$ compared to 
the true information $I_{\alpha\alpha}$.

\begin{table}[htbp]
\caption[percentage]{Observed information  (variance)\ ($n=50$, 1000 iteration)}
\begin{center}
\begin{tabular}{cccccc}    \hline
$\alpha$ & $\bar{\alpha}$  & $\hat{I}_{\alpha\alpha}$(1) & $\hat{I}_{\alpha\alpha}$(2) 
&   $I_{\alpha\alpha}$   \\ \hline
1.5  & 1.531  &  0.4863 \quad (0.033)  & 0.5145 \quad (0.018)  &  0.4737    \\
1.0  & 1.025  &  0.8541 \quad (0.157)  & 0.9001 \quad (0.099)  &  0.8590     \\
0.5  & 0.509  &  4.2819 \quad (3.517)  & 4.4660 \quad (2.755)  &
4.2748  \\ \hline
\end{tabular}
\label{tbl:observed-info-alpha}
\end{center}
\end{table}

For the rest of this section, we discuss behavior of the Fisher
information as $\alpha \uparrow 2$. Investigation of the Fisher
information in this case is difficult because it involves detailed
study of behavior of the score function as $\alpha \uparrow 2$ and
$x\rightarrow\infty$.

DuMouchel (1975, 1983) have proved that $I_{\alpha\alpha} \rightarrow
\infty$ as $\alpha \uparrow 2$.    
Nagaev  and  Shkol'nik (1988) made more detailed analysis of
$I_{\alpha\alpha}$ and shown that
\begin{equation}
\label{eq:Nagaev}
I_{\alpha\alpha} = \frac{1}{4 \Delta\log(1/\Delta)} (1+o(1)), \qquad \Delta=2-\alpha.
\end{equation}

Table \ref{tbl:information-alpha-sigma} gives the information for
$\alpha$ extremely close to 2.
$I_{\alpha\alpha}(1)$ means 
use of Taylor approximation around normal (\ref{eq:alpha-2}) 
and $I_{\alpha\alpha}(2)$
means use of (\ref{eq:xtoinf}) when $x>10^{\frac{3}{1+\alpha}}$ for the
density. The same is done with $I_{\sigma\alpha}$. There seems to be no
substantial difference between $I_{\alpha\alpha}(1)$ and  $I_{\alpha\alpha}(2)$. 
In the column N\&S we show the values of $1/(4 \Delta\log(1/\Delta))$.
Our computation for $\Delta \ge 10^{-6}$ is accurate.  The
convergence of 
(\ref{eq:Nagaev}) seems to be very slow.
Limiting theoretical behavior of $I_{\sigma\alpha}(\alpha)$ 
as $\alpha \uparrow 2$ is not known at present.

\begin{table}[htbp]
\caption[percentage]{Information around  $\alpha=2$}
\begin{center}
\begin{tabular}{llllll}    \hline
$\alpha$ & $I_{\alpha\alpha}(1)$ &  $I_{\alpha\alpha}(2)$ & N\&S & 
$I_{\sigma\alpha}(1)$ & $I_{\sigma\alpha}(2)$ 
  \\ \hline
$2.0-10^{-10}$  & 106860414 &  92384764 & 108573620 &  $-1.3482$ & $-1.3427$ \\ 
$2.0-10^{-9}$   & 10810787  &  10167389 & 12063736  &  $-1.3482$ & $-1.3395$ \\ 
$2.0-10^{-8}$   & 1144778   &  1131645    & 1357170&  $-1.3482$ & $-1.3123$  \\
$2.0-10^{-7}$   & 127953    &  127802     & 155105.2&  $-1.3478$ & $-1.2553$ \\
$2.0-10^{-6}$   & 14724     &  14722      & 18095.60&  $-1.2094$ & $-1.1923$ \\
1.99999      & 1750      &  1750       & 2171.472 &  $-1.1100$ & $-1.1100$ \\
1.9999      & 217       &  217        &  271.4341 &  $-1.0069$ & $-1.0069$ \\
\end{tabular}
\label{tbl:information-alpha-sigma}
\end{center}
\end{table}

\section{Conclusions and some discussions} 
\label{sec:discussion}

In this paper we proposed reliable numerical calculations of the
symmetric stable densities and their partial derivatives including
various boundary cases.  We found that except for very small values of
$\alpha$ ($\alpha < 0.1$) our method works very well.  This enables us
to reliably compute the maximum likelihood estimator of the symmetric stable
distributions and its standard error.

For the family of stable distributions,  the use of the
observed Fisher information in (\ref{eq:observed-FI}) for assessing the
standard deviation of the maximum likelihood estimator needs further
investigation.  Our simulation suggests that there may be some merit
in using only the first term on the right hand side of
(\ref{eq:observed-FI}).

Further study is needed to theoretically establish the limiting
behavior of the Fisher information matrix as 
$\alpha \uparrow 2$.

Finally it is of interest to extend the methods of the present paper to the
general asymmetric stable densities and to the multivariate symmetric
stable densities.  These extensions will be studied in our subsequent
works.

\appendix

\section{The third order derivative of $f(x;\alpha)$ w.r.t.\ $\alpha$
around Cauchy}
\label{sec:alpha-3}

Write 
$$
G(\nu,y)=\int_0^\infty  u^{\nu-1}e^{-u} \log^3 u \cos(uy)du.
$$
Differentiating 
(\ref{eq:oberhettinger})
three times and setting  $a=1$ we get
\begin{eqnarray}
G(\nu,y)
&=& (1+y^2)^{-\frac{1}{2}\nu}
\left\{\frac{\Gamma(\nu)}{4}\log^2(1+y^2)-\Gamma'(\nu)\log(1+y^2)+\Gamma''(\nu) \right\}
\nonumber \\
&& \hspace{1.0cm}
\times \left\{\cos(\nu z)\left(\psi(\nu)-\frac{1}{2}\log(1+y^2)\right)-z\sin(\nu)\right\} 
\nonumber\\
&& 
+(1+y^2)^{-\frac{1}{2}\nu}
\left\{-\Gamma(\nu)\log(1+y^2)+2\Gamma'(\nu)\right\} \nonumber \\
&& \hspace{1.0cm}
\times 
\left\{-z\sin(\nu z)\left(\psi(\nu)-\frac{1}{2}\log(1+y^2)\right)-z^2\cos(\nu z)+\cos(\nu z)
\psi'(\nu z) \right\} \nonumber \\
&&%
+(1+y^2)^{-\frac{1}{2}\nu}\Gamma(\nu)
\left[-z^2\cos(\nu z)\left(\psi(\nu)-\frac{1}{2}\log(1+y^2)\right) \right.\nonumber \\
&& \hspace{1.0cm}
 -2z\sin(\nu z)\psi'(\nu)
+z^3\sin(\nu z)+\cos(\nu z)\psi''(\nu)\bigg], \nonumber
\end{eqnarray}
where 
$z=\arctan y$. 
Then 
\[
 f_{\alpha\alpha\alpha}(x;1)=\frac{1}{\pi}\left(-G(4,x)+3G(3,x)-G(2,x)\right).
\]

\section{Derivatives of the density at $\alpha=0.5$}
\label{sec:alpha05}

For the purpose of checking some of our calculations, we can use the
explicit formula of the density at the special case of $\alpha=0.5$.  {}From 
(2.8.30) of Zolotarev (1986)  
$f(x;0.5)$ is written as 

\[
f(x;0.5)= \frac{x^{-\frac{3}{2}}}{\sqrt{2\pi}}\left[
\sin\left(\frac{1}{4x}\right)
\left\{\frac{1}{2}-S\left(\frac{1}{\sqrt{2\pi x}}\right)\right\}
+\cos\left(\frac{1}{4x}\right)
\left\{\frac{1}{2}-C\left(\frac{1}{\sqrt{2\pi x}}\right)\right\}
\right], 
\]
where
\[
 S(x)=\int_0^x \sin\left(\frac{\pi}{2}t^2\right) \quad \mbox{and} \quad
 C(x)=\int_0^x \cos\left(\frac{\pi}{2}t^2\right)
\]
are known as Fresnel integral functions.
We differentiate the above representation and obtain
\begin{eqnarray*}
f'(x;0.5) &=& -\frac{3}{2}\frac{x^{-\frac{5}{2}}}{\sqrt{2\pi}}\left[
\sin\left(\frac{1}{4x}\right)
\left\{\frac{1}{2}-S\left(\frac{1}{\sqrt{2\pi x}}\right)\right\}
+\cos\left(\frac{1}{4x}\right)
\left\{\frac{1}{2}-C\left(\frac{1}{\sqrt{2\pi x}}\right)\right\}\right] \\
          & & -\frac{x^{-\frac{7}{2}}}{4\sqrt{2\pi}}\left[
\cos\left(\frac{1}{4x}\right)
\left\{\frac{1}{2}-S\left(\frac{1}{\sqrt{2\pi x}}\right)\right\}
-\sin\left(\frac{1}{4x}\right)
\left\{\frac{1}{2}-C\left(\frac{1}{\sqrt{2\pi
        x}}\right)\right\}\right] 
+\frac{x^{-3}}{4\pi}
\end{eqnarray*}
and
\begin{eqnarray*}
f''(x;0.5) &=& \frac{x^{-\frac{7}{2}}}{4\sqrt{2\pi}}\left(15
-\frac{x^{-2}}{4}\right)
\left[\sin\left(\frac{1}{4x}\right)
\left\{\frac{1}{2}-S\left(\frac{1}{\sqrt{2\pi x}}\right)\right\}
+\cos\left(\frac{1}{4x}\right)
\left\{\frac{1}{2}-C\left(\frac{1}{\sqrt{2\pi x}}\right)\right\}\right] \\
          & & +\frac{5}{4}\frac{x^{-\frac{9}{2}}}{\sqrt{2\pi}}\left[
\cos\left(\frac{1}{4x}\right)
\left\{\frac{1}{2}-S\left(\frac{1}{\sqrt{2\pi x}}\right)\right\}
-\sin\left(\frac{1}{4x}\right)
\left\{\frac{1}{2}-C\left(\frac{1}{\sqrt{2\pi
        x}}\right)\right\}\right] 
-\frac{9}{8}\frac{x^{-4}}{\pi} .
\end{eqnarray*}
We have confirmed that our formulas in Section \ref{sec:derivative-location-scale}
numerically coincide  with these explicit expressions at $\alpha=1/2$ including
the boundary cases.


\begin{thebibliography}{4}

\bibitem{bergstrom} Bergstr\"om, H. (1953).
On some expansions of stable distribution functions. \textit{Ark.\
  Mat.}, \textbf{2}, 375--378.

\bibitem{Brorsen} Brorsen,\ B.\ W.\ and\ Yang,\ S.\ R.\ (1990).\ Maximum
        likelihood estimates of symmetric stable distribution
        parameters.\ \textit{Comm.\ Statist.\ Simulation Comput.},
\textbf{19},  1459--1464. 

\bibitem{dumouchel73} DuMouchel, W.\ H. (1973). 
On the asymptotic normality of the maximum-likelihood estimate when
sampling from a stable distribution. 
\textit{Ann. Statist.},  \textbf{1}, 948--957.

\bibitem{dumouchel75} DuMouchel, W.\ H. (1975). 
Stable distributions in statistical inference 2: Information from
        stably distributed samples. 
\textit{J.\ Amer.\ Statist.\ Assoc.}, \textbf{70}, 386--393.

\bibitem{dumouchel} DuMouchel, W.\ H. (1983). 
Estimating the stable index $\alpha$ in order to measure tail thickness:
        a critique.
\textit{Ann. Statist.}, \textbf{11}, 1019--1031.

\bibitem{efron-hinkley} Efron, B.\ and Hinkley, D.\ V. (1978).
Assessing the accuracy of the maximum likelihood
      estimator: Observed versus expected Fisher information.
\textit{Biometrika}, \textbf{65},
457--481.

\bibitem{Feller-vol2}
Feller, W. (1971).  \textit{An Introduction to Probability Theory and Its
  Applications.\ Vol.\ 2}, 2nd ed.,  Wiley, New York.

\bibitem{Free Software Foundation, Inc.} Free Software Foundation,\
        Inc.\ (2003).\ \textit{GNU Scientific Library - Reference
        Manual.}\ Edition 1.4,\ for GSL Version 1.4,
  59 Temple Place - Suite 330, Boston, MA 02111, USA. 

\bibitem{mcculloch} McCulloch,\ J.\ H.\ (1998).\ Numerical approximation of the
    symmetric stable distribution and density.\ \textit{A Practical Guide to
        Heavy Tails} (R.\ J. Adler et al.\ eds.),
 Birkhauser, Boston,\ 489--499.


\bibitem{nagaev-shkolnik}
Nagaev, A.\ V.\  and  Shkol'nik, S.\ M.   (1988). Some properties
  of symmetric stable distributions close to the normal
  distribution. \textit{Theory of Probability and its Applications},  
\textbf{33}, 139--144.



\bibitem{nolana} Nolan,\ J.\ P.\ (1997). Numerical calculation of 
  stable densities and distribution.\ \textit{Comm.\ Statist.\
        Stochastic Models}, \textbf{13},\ 759--774.

\bibitem{nolanb} Nolan,\ J.\ P.\ (1998) \ Parameterizations and modes
        of stable distributions.\ \textit{Statist.\ 
    Prob.\ Lett.},\ \textbf{38},\ 187-195.

\bibitem{nolanc} Nolan,\ J.\ P.\ (2001).\ Maximum likelihood estimation
        and diagnostics for stable distributions.\ \textit{L\'evy
          Processes: Theory and Applications} 
(O.\ E.\ Barndorff-Nielsen et al.\ eds.),
Birkhauser, Boston, 
  379--400.

\bibitem{oberhettinger} Oberhettinger,\ F.\ (1990).\ \textit{Tables of Fourier
    Transforms and Fourier Transforms of Distributions.}\
  Springer-Verlag, Berlin. %

\bibitem{zolotarev} Zolotarev,\ V.\ M.\ (1986).\ \textit{One-Dimensional
        Stable Distributions.} %
        Transl.\ of Math.\ Monographs,\ \textbf{65},\ Amer.\ Math.\  Soc., Providence,
        RI.\ (Transl.\ of the original 1983 Russian)



\end{thebibliography}
\end{document}